\documentclass[12pt]{article}

\usepackage {amssymb}
\usepackage {amscd}
\usepackage {latexsym}
\usepackage {euscript}

\usepackage[T2A]{fontenc}
\usepackage[cp1251]{inputenc}
\usepackage[english,russian]{babel}
\usepackage[tbtags]{amsmath}
\usepackage{amsfonts,amssymb}

\usepackage{mathrsfs}
\usepackage{graphicx}

\newtheorem{lemma}{Lemma}
\newtheorem{theorem}{Theorem}
\newtheorem{definition}{Definition}
\newtheorem{proposition}{Proposition}
\newtheorem{corollary}{Corollary}
\newtheorem{remark}{Remark}

\begin{document}

\bigskip

\bigskip

\bigskip

{\Large
\begin{center} Multipliers in the scale of periodic Bessel potential spaces with smoothness indices of different signs
 \end{center}}

\bigskip

\centerline{ALEXEI\,A.~BELYAEV, ANDREI\,A.~SHKALIKOV\footnote{This work is supported by Russian Science Foundation (RNF) under grant No 20-11-20261.}}

\bigskip

\medskip

\begin{quote}

{\large  Abstract.} We prove a general type description result for the multipliers acting between two periodic Bessel potential spaces, defined on the $n$--dimensional torus, in a case when their smoothness indices are of different sign. This is done through the detailed examination of a periodic analogue of the linear operator $J_s$, which is employed in the definition of the scale of the Bessel potential space defined on the whole space $\mathbb{R}^n$. Our method of defining this periodic analogue of $J_s$ uses the results about an asymptotic behaviour of the generalized Fourier coefficients and existence of a natural homeomorphism between the spaces $\mathcal{D}'(\mathbb{T}^n)$ and $S'_{2 \cdot \pi}(\mathbb{R}^n)$, where the latter consists of all $2 \cdot \pi$--periodic distributions from the dual Schwartz space $\mathcal{S}'(\mathbb{R}^n)$.

\end{quote}

\bigskip

\bigskip

\medskip


\bigskip

{\Large 1. Introduction.}

\bigskip

\bigskip

In this paper we consider the multipliers acting between the two periodic Bessel potential spaces, defined on the $n$--dimensional torus, in a case when the smoothness index of a first space is nonegative and the smoothness index of a second space is nonpositive. Our main aim is to obtain a natural constructive description of this multiplier space in terms of the scale of the periodic Bessel potential spaces.

Namely, our main result is the following one.

\begin{theorem}\label{periodic_multipliers_description}
Let $s \in [0; \: +\infty), \: t \in [0; \: +\infty), \: p \in (1; \; +\infty), \: q \in (1; \: +\infty)$. Additionally, 

if $s \geqslant t$, then let $s > \frac{n}{p}$ and one of the following conditions hold 
$$
1) \; p \leqslant q', \; s - \frac{n}{p} \geqslant t - \frac{n}{q'} \quad \mbox{or}  \quad 2) \; p \geqslant q',
$$

and if $t \geqslant s$, then let $t > \frac{n}{q'}$ and one of the following conditions hold 
$$
3) \; q' \leqslant p, \; t - \frac{n}{q'} \geqslant s - \frac{n}{p} \quad \mbox{or} \quad 4) \; q' \geqslant p.
$$

Then
$$
M_{\pi}[H^s_p(\mathbb{T}^n) \to H^{-t}_q(\mathbb{T}^n)] = H^{-t}_q(\mathbb{T}^n) \cap H^{-s}_{p'}(\mathbb{T}^n),
$$
and the norms $\| \cdot \|_{M_{\pi}[H^s_p(\mathbb{T}^n) \to H^{-t}_q(\mathbb{T}^n)]}$ and $\| \cdot \|_{H^{-t}_q(\mathbb{T}^n) \cap H^{-s}_{p'}(\mathbb{T}^n)}$ are equivalent. Here we define the norm $\| \cdot \|_{H^{-t}_q(\mathbb{T}^n) \cap H^{-s}_{p'}(\mathbb{T}^n)}$ by letting
$$
\forall \: u \in H^{-t}_q(\mathbb{T}^n) \cap H^{-s}_{p'}(\mathbb{T}^n) \quad \; \| u \|_{H^{-t}_q(\mathbb{T}^n) \cap H^{-s}_{p'}(\mathbb{T}^n)} = \max\left(\| u \|_{H^{-t}_q(\mathbb{T}^n)}; \: \| \cdot \|_{H^{-s}_{p'}(\mathbb{T}^n)}\right).
$$
\end{theorem}

The result of Theorem \ref{periodic_multipliers_description} can be viewed as a natural counterpart to the result about the description of multipliers in a case of the whole space $\mathbb{R}^n$, established in [2, Theorem 1], where under similar assumptions it was proved that
$$
M[H^s_p(\mathbb{T}^n) \to H^{-t}_q(\mathbb{T}^n)] = H^{-t}_{q, \: unif}(\mathbb{T}^n) \cap H^{-s}_{p', \: unif}(\mathbb{T}^n).
$$

Here the compactness of the $n$--dimensional torus allows us to use the scale of the periodic Bessel potential spaces themselves instead of the scale of the uniformly localized Bessel potential spaces $(H^{\gamma}_{r, \; unif}(\mathbb{R}^n) \: | \; \gamma \in \mathbb{R}, \: r \in (1; \; +\infty))$.

In order to establish this result we use the definition in terms of the Fourier coefficients of a natural lifting operator $J_{s, \: \pi}$, which acts on the space of all $2 \cdot \pi$--periodic distributions from the dual Schwartz space $\mathcal{S}'(\mathbb{R}^n)$. Due to the existence of a natural homeomorphism between this space and the space $\mathcal{D}'(\mathbb{T}^n)$, which is dual to the space of all infinitely differentiable functions defined on torus, it can be carefully shown that the action of this operator $J_{s, \: \pi}$ agrees well with an action of the operator $J_s \colon \mathcal{D}(\mathbb{T}^n) \to \mathcal{D}(\mathbb{T}^n)$. The fact that a function, which is identically equal to $1$ on torus, generates a regular distribution that belongs to any periodic Bessel potential space is also crucial for our proof of this constructive description.

Establishing of a constructive description of the space $M[H^s_p(\mathbb{T}^n) \to H^{-t}_q(\mathbb{T}^n)]$ is especially useful when we treat the singular perturbations of the strongly elliptic operators, acting on the scale of the periodic Bessel potential spaces, in a situation when the perturbation potential is a general type distribution from the periodic Bessel potential space with a negative index of smoothness. Studying the singular perturbations of the general strongly elliptic operators defined on the $n$--dimensional torus, among which the fractional Laplacian plays a major role, remains a heavily investigated topic (see, e.g., [3, 9, 10, 11]) which has a close connection to the problem of ill--posedness for the Navier--Stokes type equations (see [4]).

Investigating the properties of the fractional Laplacian on the $n$--dimensional torus (which corresponds to the case of periodic conditions) and its singular perturbations turns out to be very much interconnected with the embeddings and description theorems for the multiplier spaces for the periodic Bessel potential spaces. On the other hand, it allows us to get a better grasp on the specifics of this situation when compared with a general problem concerning the action of the fractional Laplacian defined on a general type smooth domain and generated by some conditions of Dirichlet and Neumann type (for the latter considerations see, e.g., [7]).

\bigskip

\bigskip

\bigskip

\bigskip

{\Large 2. Basic definitions and notation. Function spaces on the $n$--dimensional torus.}

\bigskip

\bigskip

Let $n \in \mathbb{N}$. We shall consider the $n$--dimensional torus $\mathbb{T}^n \stackrel{def}{=} \underbrace{S_1 \times \ldots \times S_1}_{n}$, where
$$
S_1 = \{ z \in \mathbb{C} \: | \; |z| = 1 \}.
$$

By $F(\mathbb{T}^n)$ we denote the set of all complex-valued functions, defined on $\mathbb{T}^n$, which becomes a linear space with respect to the pointwise operations of addition of two functions and multiplication of a function by a complex number, and $F_{2 \pi}(\mathbb{R}^n)$ shall denote the set of all $2 \pi$--periodic (in each of its $n$ variables) functions, defined on $\mathbb{R}^n$, which is also a linear space with respect to the pointwise linear operations.

The mapping
$$
\Pi \colon \mathbb{R}^n \to \mathbb{T}^n; \; \; (x_1, \ldots, x_n) \stackrel{\Pi}{\longmapsto} (e^{i x_1}, \ldots, e^{i x_n})
$$
naturally induces the mapping $\mathbf{\Pi} \colon F(\mathbb{T}^n) \to F_{2 \pi}(\mathbb{R}^n)$, defined as follows:

for any function $f \in F(\mathbb{T}^n)$ we define a function $\mathbf{\Pi}(f) \colon \mathbb{R}^n \to \mathbb{C}$ (denoted below by $f_{\pi}$), where
$$
\forall \: x \in \mathbb{R}^n \; \; f_{\pi}(x) = f(\Pi(x)).
$$

It is easy to see that $Im_{\mathbf{\Pi}}(F(\mathbb{T}^n)) = F_{2 \pi}(\mathbb{R}^n)$.

In what follows for any $\alpha = (\alpha_1, \ldots, \alpha_n) \in \mathbb{Z}^n$ we use the following notations:
$$
| \alpha |_1 = \sum\limits_{k = 1}^n |\alpha_k|, \; \; \; | \alpha | = \biggl( \sum\limits_{k = 1}^n \left|\alpha_k\right|^2 \biggr)^{\frac{1}{2}}.
$$

We say that a function $f \colon \mathbb{T}^n \to \mathbb{C}$ is differentiable at $z^{(0)} \in \mathbb{T}^n$, if the function $f_{\pi}$ is differentiable at an arbitrary point $x \in \mathbb{R}^n$, such that $\Pi(x) = z^{(0)}$, and for arbitrary $k \in \overline{1, \: n}$ operator of the periodic partial differentiation on $\mathbb{T}^n$ is defined by
$$
\partial_k f(z) \stackrel{def}{=} \left(\frac{\partial}{\partial x_k} f_{\pi}\right)(x),
$$
where $z \in \mathbb{T}^n$ and $x \:$ is an arbitrary element of $\mathbb{R}^n$, such that $\Pi(x) = z$.

For any multi-index $\alpha = (\alpha_1, \ldots, \alpha_n) \in \mathbb{Z}_+^{\: n}$ a periodic partial differential operator $D^{\alpha}$ is defined in exactly the same way for such functions $f\in F(\mathbb{T}^n)$ that the action of the partial differential operator $\frac{\partial^{\alpha_1}}{(\partial x_1)^{\alpha_1}} \ldots \frac{\partial^{\alpha_n}}{(\partial x_n)^{\alpha_n}}$ is well--defined on $f_{\pi}$.

The set $\mathcal{D}(\mathbb{T}^n)$ of all complex--valued infinitely differentiable functions, defined on the $n$--dimensional torus $\mathbb{T}^n$, which becomes a linear space with respect to the pointwise linear operations, is endowed with a countable family of norms $\{\| \cdot \|_m \}_{m \in \mathbb{Z}_+}$, where for arbitrary $m \in \mathbb{Z}_+$ the norm $\| \cdot \|_m$ is defined by
$$
\forall \; f \in D(\mathbb{T}^n) \; \; \: \| f \|_m \stackrel{def}{=} \max\limits_{\substack{\alpha \in \mathbb{Z}_+^{\: n} \colon \\ |\alpha| \leqslant m}} \left(\max_{z \in \mathbb{T}^n} \left|(D^{\alpha}f)(z)\right|\right).
$$
Let us remark that an image of $\mathcal{D}(\mathbb{T}^n)$ under the mapping
$$
\Pi \colon F(\mathbb{T}^n) \to F_{2 \pi}(\mathbb{R}^n), \; \; f \stackrel{\Pi}{\longmapsto} f_{\pi},
$$
coincides with the set $C^{\infty}_{2 \pi}(\mathbb{R}^n)$ of all infinitely differentiable complex--valued $2 \pi$--periodic functions defined on $\mathbb{R}^n$.

Since the family of norms $\{\| \cdot \|_m \}_{m \in \mathbb{Z}_+}$, which induce the natural polynormed topology $\tau_{\mathcal{D}(\mathbb{T}^n)}$ on $\mathcal{D}(\mathbb{T}^n)$, is countable, this topology is metrizable and, hence, for any topological space $(X, \: \tau)$ the notions of continuity and sequential continuity for the mappings from $\mathcal{D}(\mathbb{T}^n)$ to $X$ (with respect to the topology, generated by this countable family of norms, and the topology $\tau$) are equivalent. Therefore, a linear mapping $u \colon \mathcal{D}(\mathbb{T}^n) \to \mathbb{C}$ belongs to the dual space $\mathcal{D}'(\mathbb{T}^n)$, which consists of all complex--valued continuous (with respect to the topology $\tau_{\mathcal{D}(\mathbb{T}^n)}$) linear functionals from $\mathcal{D}(\mathbb{T}^n)$ to $\mathbb{C}$, if and only if for any $f_0 \in \mathcal{D}(\mathbb{T}^n)$ and for any sequence $(f_n \in \mathcal{D}(\mathbb{T}^n) \: | \; n \in \mathbb{N})$, such that $f_n \xrightarrow[n \to \infty]{\tau_{\mathcal{D}(\mathbb{T}^n)}} f_0$, we have $\lim\limits_{n \to \infty} (u(f_n)) = u(f_0)$.

This dual linear space $\mathcal{D}'(\mathbb{T}^n)$ is endowed with the weak-$*$ topology, denoted by $\tau_{\mathcal{D}'(\mathbb{T}^n)}$, which is generated by the family of seminorms $\{\| \cdot \|_f \}_{f \in D(\mathbb{T}^n)}$, where for any function $f \in D(\mathbb{T}^n)$ seminorm $\| \cdot \|_f$ is defined as follows:
$$
\forall \: u \in D'(\mathbb{T}^n) \; \: \| u \|_f = |u(f)|.
$$
Thus, a sequence of distributions $(u_k \in D'(\mathbb{T}^n) \: | \; k \in \mathbb{N})$ converges to a distribution $u_0 \in D'(\mathbb{T}^n)$ (which will be denoted as $u_k \xrightarrow[k \to \infty]{D'(\mathbb{T}^n)} u$ in what follows) if and only if
$$
\forall \: f \in D(\mathbb{T}^n) \; \; \: u_k(f) \xrightarrow[k \to \infty]{} u(f).
$$

For an arbitrary number $p \in [1; \: +\infty)$ we say that a function $f \colon \mathbb{T}^n \to \mathbb{C}$ belongs to $\mathcal{L}_p(\mathbb{T}^n)$, if the function $|f_{\pi}|^p$ is integrable with respect to a restriction of the classical Lebesgue measure $\mu_{cl}$ onto $[-\pi; \: \pi]^n$. The set $\mathcal{L}_p(\mathbb{T}^n)$ is a complex linear space with respect to the pointwise linear operations and the function $\| \cdot \|_{\mathcal{L}_p(\mathbb{T}^n)} \colon \mathcal{L}_p(\mathbb{T}^n) \to \mathbb{R}$, defined by
$$
\forall \; f \in \mathcal{L}_p(\mathbb{T}^n) \; \; \| f \|_{\mathcal{L}_p(\mathbb{T}^n)} = \left(\int\limits_{[-\pi; \: \pi]^n} |f_{\pi}(x)|^p \: d\mu_{cl}(x)\right)^{\frac{1}{p}},
$$
is a seminorm on this linear space. Employing a general construction of the quotient space of a seminormed space by the kernel of its seminorm, we obtain the normed space $(L_p(\mathbb{T}^n), \: \| \cdot \|_{L_p(\mathbb{T}^n)})$, where $L_p(\mathbb{T}^n) = \left\{ [f] \: | \; f \in \mathcal{L}_p(\mathbb{T}^n) \right\}$ and for any $f \in \mathcal{L}_p(\mathbb{T}^n)$ the equivalence class $[f]$ is defined by
$$
[f] = \left\{ g \in \mathcal{L}_p(\mathbb{T}^n) \: | \; \| g - f \|_{\mathcal{L}_p(\mathbb{T}^n)} = 0 \right\} =
$$
$$
= \left\{ g \in \mathcal{L}_p(\mathbb{T}^n) \: | \; \mu_{L, \: cl}(\{x \in [-\pi; \: \pi]^n \: | \; g_{\pi}(x) \neq f_{\pi}(x) \}) = 0 \right\},
$$
linear operations on $L_p(\mathbb{T}^n)$ are defined naturally by applying the linear operations on arbitrary representative functions of the corresponding equivalence classes and the norm $\| \cdot \|_{L_p(\mathbb{T}^n)} \colon L_p(\mathbb{T}^n) \to \mathbb{R}$ is defined as follows
$$
\forall \: f \in \mathcal{L}_p(\mathbb{T}^n) \; \; \| [f] \|_{L_p(\mathbb{T}^n)} = \| f \|_{\mathcal{L}_p(\mathbb{T}^n)}.
$$

For an arbitrary number $p \in [1; \: +\infty)$ in the sequel we shall often identify a function $f \in \mathcal{L}_p(\mathbb{T}^n)$ and an equivalence class $[f] \in L_p(\mathbb{T}^n)$, generated by $f$.

Now let us define the family of functions $\Phi = \{ f_k \: | \; k \in \mathbb{Z}^n \} \subset \mathcal{D}(\mathbb{T}^n)$. For an arbitrary $k = (k_m \in \mathbb{Z} \: | \; m \in \overline{1, \: n}) \in \mathbb{Z}^n$ we define the function $f_k \colon \mathbb{T}^N \to \mathbb{C}$ by letting
$$
\mbox{for any} \; z = (z_m \in S_1 \: | \; m \in \overline{1, \: n}\,) \in \mathbb{T}^n \; \; f_k(z) = \frac{z^k}{(2 \cdot \pi)^{\frac{n}{2}}} \: , \; \; \mbox{where} \; \; z^k \: \stackrel{def}{=} \prod\limits_{m = 1}^n z_m^{k_m} \: .
$$
Since for any $k \in \mathbb{Z}^n$ for the function $g_k \in F_{2 \cdot \pi}(\mathbb{R}^n)$, defined by
$$
\forall \: x \in \mathbb{R}^n \; \; \; g_k(x) \stackrel{def}{=} \frac{1}{(2 \cdot \pi)^{\frac{n}{2}}} \: e^{i \cdot \: < k, x >},
$$
we have $\mathbf{\Pi}(g_k) = f_k$, then, as it is well--known (see, e.g., [8, Theorem B.1]), the family $\Phi$ (or, more precisely, the family of equivalence classes $\{[f_k] \: | \; k \in \mathbb{Z}^n\}$) forms a countable orthonormal basis in the Hilbert space $\left(L_2(\mathbb{T}^n), \: < \cdot \, ; \: \cdot >_{L_2(\mathbb{T}^n)}\right)$.

Here we outline the construction of a natural embedding of $\mathcal{L}_1(\mathbb{T}^n)$ into the dual space $\mathcal{D}'(\mathbb{T}^n)$, which plays a crucial role in what follows. For any function $f \in \mathcal{L}_1(\mathbb{T}^n)$ we define a functional $\mathbf{f} \colon \mathcal{D}(\mathbb{T}^n) \to \mathbb{C}$ by letting
$$
\forall \: \varphi \in D(\mathbb{T}^n) \; \; \; \: \mathbf{f}(\varphi) = \int\limits_{[-\pi; \: \pi]^n} f_{\pi}(x) \cdot \varphi_{\pi}(x) \: d\mu_{cl}(x).
$$
It is well--known that this functional $\mathbf{f}$, generated by $f \in \mathcal{L}_1(\mathbb{T}^n)$, is well--defined as a regular distribution belonging to the dual space $\mathcal{D}'(\mathbb{T}^n)$ and that for any $f_1 \in [f]$ regular distributions $\mathbf{f}$ and $\mathbf{f_1}$, so we can define a regular distribution, generated by the equivalence class $[f] \in L_1(\mathbb{T}^n)$.

Since for any $p \in [1; \: +\infty)$ we have the embedding $\mathcal{L}_p(\mathbb{T}^n) \subset \mathcal{L}_1(\mathbb{T}^n)$ (which is continuous with respect to the topologies, generated by the seminorms $\| \cdot \|_{\mathcal{L}_p(\mathbb{T}^n)}$ and $\| \cdot \|_{\mathcal{L}_1(\mathbb{T}^n)}$ respectively), then for any function $f \in \mathcal{L}_p(\mathbb{T}^n)$ (where $p \in [1; \: +\infty)$) the regular distribution $\mathbf{f} \in \mathcal{D}'(\mathbb{T}^n)$ is well--defined (and this regular distribution remains the same if we construct it from some other function from the equivalence class $[f]$), which allows us to guarantee the validity of a natural embedding of both $\mathcal{L}_p(\mathbb{T}^n)$ and $L_p(\mathbb{T}^n)$ into $\mathcal{D}'(\mathbb{T}^n)$ by identifying the function $f \in \mathcal{L}_p(\mathbb{T}^n)$ (or the equivalence class $[f] \in L_p(\mathbb{T}^n)$) with the regular distribution $\mathbf{f} \in \mathcal{D}'(\mathbb{T}^n)$.

By $\mathbf{D}(\mathbb{T}^n)$ we shall denote the set of all regular distributions from $\mathcal{D}'(\mathbb{T}^n)$ with their generating functions being elements of $\mathcal{D}(\mathbb{T}^n)$.

\bigskip

\bigskip

\bigskip

\bigskip

{\Large 3. Classical facts on the periodic Fourier decompositions.}

\bigskip

\bigskip

Let us present some classical facts from harmonic analysis which will be needed in the sequel.

Let $k \in \mathbb{Z}^n$. Then for an arbitrary function $f \in \mathcal{D}(\mathbb{T}^n)$ we shall denote by $C_k(f)$ its $k$--th Fourier coefficient with respect to the orthonormal (with respect to the natural positive semi--definite inner product, defined on $\mathcal{L}_2(\mathbb{T}^n)$) system $\Phi$, defined as
$$
C_k(f) = \int\limits_{[-\pi; \: \pi]^n} f_{\pi}(x) \cdot \overline{g_k(x)} \: d\mu_{cl}(x).
$$
On the other hand, for an arbitrary distribution $u \in \mathcal{D}'(\mathbb{T}^n)$ we shall use the notation $c_k(u) \stackrel{def}{=} u(f_{-k})$ and call the number $c_k(u) \in \mathbb{C}$ as a $k$--th distributional Fourier coefficient of the distribution $u$.

There are two classical results on the convergence of Fourier series, defined by these two sequences of Fourier coefficients.

\begin{proposition}\label{periodic_function_Fourier_coefficients}(see, e.g., [14, 3.2.2.(i)])

Let $f \in \mathcal{D}(\mathbb{T}^n)$.
Then the series $\sum\limits_{k \in \mathbb{Z}^n} C_k(f) \cdot f_k$ converges to the function $f$  with respect to the topology $\tau_{\mathcal{D}(\mathbb{T}^n)}$.
\end{proposition}

\begin{proposition}\label{periodic_distribution_Fourier_coefficients}(see, e.g., [14, 3.2.2.(ii)])

Let $u \in \mathcal{D}'(\mathbb{T}^n)$.
Then the series $\sum\limits_{k \in \mathbb{Z}^n} c_k(u) \cdot \mathbf{f}_k$ converges to $u$ (with respect to $\tau_{\mathcal{D}'(\mathbb{T}^n)}$), 
where the sequence $(c_k(u) \in \mathbb{C} \: | \; k \in \mathbb{Z}^n)$ has at most polynomial growth, which means that there exists a natural number $m \in \mathbb{N}$, such that for some $K_m > 0$ we have
$$
\forall \: k \in \mathbb{Z}^n \; \; \: |c_k(u)| \leqslant K_m \cdot |k|^m.
$$
The converse statement is also valid: namely, for any sequence $(d_k \in \mathbb{C} \: | \; k \in \mathbb{Z}^n)$ of at most polynomial growth the sequence of distributions
$$
\left(\left.\sum\limits_{\substack{k \in \mathbb{Z}^n \colon \\ |k| \leqslant m}} d_k \cdot \mathbf{f}_k \in \mathcal{D}'(\mathbb{T}^n) \: \right| \; m \in \mathbb{N}\right)
$$
converges (with respect to $\tau_{\mathcal{D}'(\mathbb{T}^n)}$) to some distribution $u_0 \in \mathcal{D}'(\mathbb{T}^n)$  and, moreover, for arbitrary multi-index $k \in \mathbb{Z}^n$ we have $c_k(u_0) = d_k$.
\end{proposition}

\begin{remark}\label{D'_Schauder_basis}
Let us note that Proposition \ref{periodic_distribution_Fourier_coefficients} implies that the family of distributions $\left(\mathbf{f}_k \in \mathcal{D}'(\mathbb{T}^n) \: | \; k \in \mathbb{Z}^n \right)$ constitutes a Schauder basis in the locally convex linear space $\left(\mathcal{D}'(\mathbb{T}^n), +, \: \cdot\right)$, whose topology is generated by the family of seminorms $\left\{\| \cdot \|_{f} \right\}_{f \in \mathcal{D}(\mathbb{T}^n)}$.
\end{remark}

\begin{remark}\label{relation_between_Fourier_coefficients}
Let us also note that Propositions \ref{periodic_function_Fourier_coefficients} and \ref{periodic_distribution_Fourier_coefficients} agree with each another.
Namely, since for an arbitrary function $f \in \mathcal{L}_2(\mathbb{T}^n)$ we have the following chain of equalities
$$
C_k(f) 
\: = \int\limits_{[-\pi; \: \pi]^n} f_{\pi}(x) \cdot \overline{e^{i \cdot <k \,; \: x>}} \: d\mu_{cl}(x) \: =
$$
$$
= \int\limits_{[-\pi; \: \pi]^n} f_{\pi}(x) \cdot e^{-i \cdot <k \,; \: x>} \: d\mu_{cl}(x) \: = \mathbf{f}(f_{-k}) \stackrel{def}{=} c_k(\mathbf{f}),
$$
it follows that for an arbitrary function $f \in \mathcal{D}(\mathbb{T}^n)$ the Fourier decompositions of $f$ itself (in $\mathcal{D}(\mathbb{T}^n)$) and of the regular distribution $\mathbf{f}$ (in $\mathcal{D}'(\mathbb{T}^n)$) have the same coefficients:
$$
f \stackrel{\mathcal{D}(\mathbb{T}^n)}{=} \sum\limits_{k \in \mathbb{Z}^n} C_k(f) \cdot f_k \quad \mbox{and} \quad \; \: \mathbf{f} \stackrel{\mathcal{D}'(\mathbb{T}^n)}{=} \sum\limits_{k \in \mathbb{Z}^n} C_k(f) \cdot \mathbf{f}_k \, .
$$
\end{remark}

\begin{remark}\label{Fourier_coefficients_fast_convergence}
There is a well--known relation between the asymptotic behaviour of the Fourier coefficients of a function $f \in \mathcal{L}_1(\mathbb{T}^n)$ and the belonging of this function $f$ to the classical (periodic) function spaces (see, e.g., [5, 3.3.2, Table 3.1]). Specifically, we need the fact that if the function $f \colon \mathbb{T}^n \to \mathbb{C}$ is integrable on $\mathbb{T}^n$, then:

$a)$ $f \in \mathcal{L}_2(\mathbb{T}^n)$ if and only if $(C_k(f) \: | \; k \in \mathbb{Z}^n) \in l_2(\mathbb{C})$;

$b)$ $f \in \mathcal{D}(\mathbb{T}^n)$ if and only if $(C_k(f) \: | \; k \in \mathbb{Z}^n)$ is a subpolynomial fast decreasing sequence, which means that
$$
\forall \: N \in \mathbb{N} \; \; \exists \: K_N > 0: \; \; \forall \: k \in \mathbb{Z}^n \: \; |C_k(f)| \leqslant \frac{K_N}{1 + |k|^N}.
$$
\end{remark}

\begin{remark}\label{Fourier_coefficients_inverse_fast_convergence}
Remark \ref{Fourier_coefficients_fast_convergence} admits a generalization in the following sense (see [14, 3.2.2]). Let $(C_k \in \mathbb{С} \: | \; k \in \mathbb{Z}^n)$ be a complex-valued sequence. Then:

$1)$ if $(C_k \in \mathbb{C} \: | \; k \in \mathbb{Z}^n) \in l_2(\mathbb{C})$, then there exists a function $f \in \mathcal{L}_2(\mathbb{T}^n)$ such that the series $\sum\limits_{k \in \mathbb{Z}^n} C_k \cdot f_k$ converges to $f$ (in a sense of the natural convergence in $\mathcal{L}_2(\mathbb{T}^n)$, generated by the seminorm $\| \cdot \|_{\mathcal{L}_2(\mathbb{T}^n)}$);

$2)$ if a sequence $(C_k \in \mathbb{C} \: | \; k \in \mathbb{Z}^n)$ decreases faster than any polynomial, then there exists a function $f \in \mathcal{D}(\mathbb{T}^n)$, such that the series $\sum\limits_{k \in \mathbb{Z}^n} C_k \cdot f_k$ converges to $f$ (in terms of the natural convergence in $\mathcal{L}_2(\mathbb{T}^n)$).

Moreover, in both $1)$ and $2)$ for any $k \in \mathbb{Z}^n$ we have
$$
C_k = \frac{1}{(2 \pi)^{\frac{n}{2}}} \cdot\int\limits_{[-\pi; \: \pi]^n} f_{\pi}(x) \cdot g_k(x) \: d\mu_{cl}(x) \stackrel{def}{=} C_k(f).
$$
\end{remark}

\bigskip

\bigskip

\bigskip

\bigskip

{\Large 4. Some additional facts on the periodic distributions and the operator $J_{s, \: \pi}$.}

\bigskip

\bigskip

For the sake of the exposition's comprehensiveness we present a proof of the following simple technical fact about the action of a distribution on $f \in \mathcal{D}(\mathbb{T}^n)$.

\begin{lemma}\label{action_of_distribution_on_function}
Let $u \in \mathcal{D}'(\mathbb{T}^n)$ and $f \in \mathcal{D}(\mathbb{T}^n)$. Then
$$
u(f) = \sum\limits_{k \in \mathbb{Z}^n} c_k(u) \cdot C_{-k}(f),
$$
where, as above, for any multi--index $k \in \mathbb{Z}^n$ we let
$$
c_k(u) = u(f_{-k}) \quad \mbox{and} \quad C_k(f) = \int\limits_{[-\pi \, ; \: \pi]^n} f_{\pi}(x) \cdot g_k(x) \: d\mu_{cl}(x).
$$
\end{lemma}

Proof. Indeed, since the distribution $u$ is sequentially continuous on $\mathcal{D}(\mathbb{T}^n)$ and a weak--$*$ convergence of the sequence $(u_k \in \mathcal{D}'(\mathbb{T}^n) \: | \; k \in \mathbb{Z}^n)$ to $u_0 \in \mathcal{D}'(\mathbb{T}^n)$ means that for an arbitrary function $\varphi \in \mathcal{D}(\mathbb{T}^n)$ the complex-valued sequence $(u_k(\varphi) \: | \; k \in \mathbb{Z}^n)$ should converge to $u_0(\varphi)$, we obtain the chain of equalities
$$
u(f) = u\left((\mathcal{D}(\mathbb{T}^n))-\sum\limits_{k \in \mathbb{Z}^n} C_k(f) \cdot f_k\right) = \sum\limits_{k \in \mathbb{Z}^n} C_k(f) \cdot u(f_k) =
$$
$$
= \sum\limits_{k \in \mathbb{Z}^n} C_k(f) \cdot \left((\mathcal{D}'(\mathbb{T}^n))-\sum\limits_{l \in \mathbb{Z}^n} c_l(u) \cdot \mathbf{f}_l\right)(f_k) =  \sum\limits_{k \in \mathbb{Z}^n} C_k(f) \cdot \left(\sum\limits_{l \in \mathbb{Z}^n} c_l(u) \cdot \mathbf{f}_l(f_k)\right).
$$
Since for any multi--index $k \in \mathbb{Z}^n$ и $l \in \mathbb{Z}^n$ we have
$$
\mathbf{f}_l(f_k) = \int\limits_{[-\pi \, ; \: \pi]^n} (f_l)_{\pi}(x) \cdot (f_k)_{\pi}(x) \: d\mu_{cl}(x) =
$$
$$
= \int\limits_{[-\pi \, ; \: \pi]^n} g_l(x) \cdot g_k(x) \: d\mu_{cl}(x) = \int\limits_{[-\pi \, ; \: \pi]^n} \frac{1}{(2 \pi)^{\frac{n}{2}}} \cdot e^{i \cdot < k \, ; \: x >} \cdot \frac{1}{(2 \pi)^{\frac{n}{2}}} \cdot e^{i \cdot < l \, ; \: x >} \: d\mu_{cl}(x) =
$$
$$
= \frac{1}{(2 \pi)^n} \cdot \int\limits_{[-\pi \, ; \: \pi]^n} e^{i \cdot < k + l \, ; \: x >} \: d\mu_{cl}(x) = \begin{cases}
1, \; \: \mbox{if} \; l = -k\\
0, \; \: \mbox{if} \; l \neq -k
\end{cases} \; ,
$$
then
$$
u(f) = \sum\limits_{k \in \mathbb{Z}^n} C_k(f) \cdot c_{-k}(u) = \sum\limits_{k \in \mathbb{Z}^n} c_k(u) \cdot C_{-k}(f).
$$

This concludes the proof of Lemma \ref{action_of_distribution_on_function}.

\bigskip

Let us recall a definition of the shift operator.

\begin{definition}\label{shift_operator}
Let $h \in \mathbb{R}^n$. Then

$1)$ the (functional) shift operator $\tau_h \colon F(\mathbb{R}^n) \to F(\mathbb{R}^n)$ is defined as follows: 

for an arbitrary function $f \in F(\mathbb{R}^n)$ we let
$$
\forall \: x \in \mathbb{R}^n \; \; (\tau_h(f))(x) = f(x - h);
$$

$2)$ the (distributional) shift operator $T_h \colon D'(\mathbb{R}^n) \to D'(\mathbb{R}^n)$ is defined as follows: 

for an arbitrary distribution $u \in D'(\mathbb{R}^n)$ we let
$$
\forall \: f \in D(\mathbb{R}^n) \; \; (T_h(u))(f) = u(\tau_{-h}(f)).
$$
\end{definition}

In what follows for an arbitrary $h \in \mathbb{R}^n$, an arbitrary function $f \in F(\mathbb{R}^n)$ and an arbitrary distribution $u \in \mathcal{D}'(\mathbb{R}^n)$ we shall often denote $\tau_h(f)$ by $f_{(h)}$ and $T_h(u)$ by $u_{(h)}$.

Let us remind that a distribution $u \in \mathcal{D}'(\mathbb{R}^n)$ is called $2 \cdot \pi$--periodic, if for any multi--index $k \in \mathbb{Z}^n$ the distributions $u_{(2 \cdot \pi \cdot k)}$ and $u$ are equal.

In the sequel 
by $\mathcal{S}^{'}_{2 \cdot \pi}(\mathbb{R}^n)$ we shall denote the set of all $2 \cdot \pi$--periodic distributions from $\mathcal{S}'(\mathbb{R}^n)$.

\begin{proposition}\label{periodization_of_toric_D_prime}(see, e.g., [14, 3.2.2(ii) and Proposition 3.2.3])

Let $u \in \mathcal{S}'(\mathbb{R}^n)$. Then the following conditions are equivalent:
$$
1) \; u \in \mathcal{S}'_{2 \pi}(\mathbb{R}^n);
$$
$$
2) \; \mbox{there exists a sequence} \; \: (d_k \in \mathbb{C} \: | \; k \in \mathbb{Z}^n) \; \: \mbox{of at most polynomial growth, such that}
$$
$$
u \stackrel{\mathcal{S}'(\mathbb{R}^n)}{=} \sum\limits_{k \in \mathbb{Z}^n} d_k \cdot \mathbf{g_k} \: ;
$$
$$
3) \; \; \mbox{there exists a distribution} \; v \in \mathcal{D}'(\mathbb{T}^n), \; \: \mbox{such that a validity of the equality}
$$
$$
v \stackrel{\mathcal{D}'(\mathbb{T}^n)}{=} \sum\limits_{k \in \mathbb{Z}^n} D_k \cdot \mathbf{f}_k
$$
$$
\mbox{for some sequence} \; \: (D_k \in \mathbb{C} \: | \; k \in \mathbb{Z}^n) \; \; \mbox{implies a validity of the equality}
$$
$$
u \stackrel{\mathcal{S}'(\mathbb{R}^n)}{=} \sum\limits_{k \in \mathbb{Z}^n} D_k \cdot \mathbf{g}_k.
$$
\end{proposition}

Moreover, it is easy to see that the statement, analogous to the converse statement of Proposition \ref{periodic_distribution_Fourier_coefficients}, also holds true:

\begin{corollary}\label{periodic_Schwartz_distribution_Fourier_coefficients}
Let a sequence $(d_k \in \mathbb{C} \: | \; k \in \mathbb{Z}^n)$ have at most polynomial growth. Then there exists a periodic distribution $u_0 \in \mathcal{S}'_{2 \pi}(\mathbb{R}^n)$ such that the series $\sum\limits_{k \in \mathbb{Z}^n} d_k \cdot \mathbf{g}_k$ converges to $u_0$ in $\mathcal{S}'(\mathbb{R}^n)$ (with respect to the weak--$*$ topology, defined on $\mathcal{S}'(\mathbb{R}^n)$).
\end{corollary}

From Proposition \ref{periodic_distribution_Fourier_coefficients} and Corollary \ref{periodic_Schwartz_distribution_Fourier_coefficients} immediately follows the existence of a natural homeomorphism $G$ between the spaces $\mathcal{D}'(\mathbb{T}^n)$ and $\mathcal{S}'_{2 \pi}(\mathbb{R}^n)$ (with respect to their natural topologies). Moreover, this homeomorphism $G \colon \mathcal{D}'(\mathbb{T}^n) \to \mathcal{S}'_{2 \pi}(\mathbb{T}^n)$ can be explicitly described in terms of the (generalized) Fourier coefficients.

Indeed, let us fix an arbitrary distribution $u \in D'(\mathbb{T}^n)$. By Proposition \ref{periodic_distribution_Fourier_coefficients}, we have
$$
u \stackrel{\mathcal{D}'(\mathbb{T}^n)}{=} \sum\limits_{k \in \mathbb{Z}^n} c_k(u) \cdot \mathbf{f}_k \: ,
$$
and the sequence $(c_k(u) \in \mathbb{C} \: | \; k \in \mathbb{Z}^n)$ has at most polynomial growth. Then, by Corollary \ref{periodic_Schwartz_distribution_Fourier_coefficients}, there exists a periodic distribution $\widetilde{u} \in \mathcal{S}'_{2 \pi}(\mathbb{R}^n)$, such that the series $\sum\limits_{k \in \mathbb{Z}^n} c_k(u) \cdot \mathbf{g}_k$ converges to $\widetilde{u}$ in the weak--$*$ topology of $\mathcal{S}'(\mathbb{R}^n)$. Then we let $G(u) = \widetilde{u}$. This equality can also be rewritten in the following manner:
$$
G(u) \stackrel{\mathcal{S}'(\mathbb{R}^n)}{=} \sum\limits_{k \in \mathbb{Z}^n} c_k(u) \cdot \mathbf{g}_k \: .
$$

\medskip

Let us now recall the definition of the operator $J_s$ in  the dual Schwartz space $\mathcal{S}'(\mathbb{R}^n)$. First of all, the Fourier transform $\mathcal{F} \colon \mathcal{S}(\mathbb{R}^n) \to \mathcal{S}(\mathbb{R}^n)$ and the inverse Fourier transform $\mathcal{F}^{-1} \colon \mathcal{S}(\mathbb{R}^n) \to \mathcal{S}(\mathbb{R}^n)$ are defined in the following way: for an arbitrary function $f \in \mathcal{S}(\mathbb{R}^n)$ we let
$$
\forall \: x \in \mathbb{R}^n \; \; (\mathcal{F}(f))(x) = \frac{1}{(2 \pi)^{\frac{n}{2}}} \cdot \int\limits_{\mathbb{R}^n} e^{-i \cdot < y \, ; \: x >} \cdot f(y) \: d\mu_{cl}(y)
$$
\centerline{and}
$$
\forall \: x \in \mathbb{R}^n \; \; (\mathcal{F}^{-1}(f))(x) = \frac{1}{(2 \pi)^{\frac{n}{2}}} \cdot \int\limits_{\mathbb{R}^n} e^{i \cdot < y \, ; \: x >} \cdot f(y) \: d\mu_{cl}(y),
$$

First of all, a classical result from harmonic analysis (see, e.g., [5, Corollary 2.2.15]) states that both Fourier transform $\mathcal{F}$ and $\mathcal{F}^{-1}$ are homeomorphisms of the Schwartz space $\mathcal{S}(\mathbb{R}^n)$ onto itself (with respect to its natural topology, induced by a countable family of seminorms).

This fact implies that the distributional analogues of the Fourier transform and the inverse Fourier transform, i.e. the linear mappings $\mathcal{F} \colon \mathcal{S}'(\mathbb{R}^n) \to \mathcal{S}'(\mathbb{R}^n)$ and $\mathcal{F}^{-1} \colon \mathcal{S}'(\mathbb{R}^n) \to \mathcal{S}'(\mathbb{R}^n)$, defined by
$$
\forall \; \varphi \in S(\mathbb{R}^n) \; \; (\mathcal{F}(u))(\varphi) \stackrel{def}{=} u(\mathcal{F}(\varphi)) \; \; \mbox{and} \; \; (\mathcal{F}^{-1}(u))(\varphi) \stackrel{def}{=} u(\mathcal{F}^{-1}(\varphi)),
$$
are also homeomorphisms on the dual Schwartz space $\mathcal{S}'(\mathbb{R}^n)$ (with respect to its weak--$*$ topology), which is a fundamental result in the distribution theory (see an analogous result in [6, Theorem 5.17]).

Then let us fix an arbitrary infinitely differentiable function $f \colon \mathbb{R}^n \to \mathbb{C}$, such that this function and any of its derivatives have at most polynomial growth. We shall use an easily verifiable fact that $f$ induces a linear operator $\mathbb{A}_f \colon \mathcal{S}(\mathbb{R}^n) \to \mathcal{S}(\mathbb{R}^n)$ of the pointwise multiplication by $f$, defined for arbitrary function $\varphi \in \mathcal{S}(\mathbb{R}^n)$ by
$$
\forall \: x \in \mathbb{R}^n \; \; (\mathbb{A}_f(\varphi))(x) = f(x) \cdot \varphi(x),
$$
and that this operator $\mathbb{A}_f$ is continuous with respect to the natural topology of $\mathcal{S}(\mathbb{R}^n)$. From this fact it easily follows that the mapping $A_f \colon \mathcal{S}'(\mathbb{T}^n) \to \mathcal{S}'(\mathbb{T}^n)$, where for arbitrary $u \in \mathcal{S}'(\mathbb{R}^n)$ we define
$$
\forall \: \varphi \in \mathcal{S}(\mathbb{R}^n) \; \; \: (A_{f}(u))(\varphi) \stackrel{def}{=} u(\mathbb{A}_f(\varphi)) = u (f \cdot \varphi),
$$
is well--defined as a linear operator on $\mathcal{S}'(\mathbb{R}^n)$ and, moreover, this linear operator $A_f$ is a homeomorphism of the dual Schwartz space $\mathcal{S}'(\mathbb{R}^n)$ onto itself (with respect to the natural weak--$*$ topology on $\mathcal{S}'(\mathbb{R}^n)$).

\begin{definition}\label{J_s_on_dual_Schwartz}
Let $s \in \mathbb{R}$. Let us also define the function $\varphi_s \colon \mathbb{R}^n \to \mathbb{C}$ as follows:
$$
\forall \: x \in \mathbb{R}^n \; \; \: \varphi_s(x) \stackrel{def}{=} \left(1 + \left|x\right|^2\right)^{\frac{s}{2}}.
$$
Then the linear operator $\mathbb{J}_s \colon \mathcal{S}(\mathbb{R}^n) \to \mathcal{S}(\mathbb{R}^n)$ is defined by
$$
\forall \: f \in S(\mathbb{R}^n) \; \; \: \mathbb{J}_s(f) \stackrel{def}{=} (\mathcal{F}^{-1} \circ \mathbb{A}_{\varphi_s} \circ \mathcal{F})(f) = \mathcal{F}^{-1}(\varphi_s \cdot \mathcal{F}(f)).
$$
and the linear operator $J_s \colon \mathcal{S}'(\mathbb{R}^n) \to \mathcal{S}'(\mathbb{R}^n)$ is defined by
$$
\forall \: u \in S'(\mathbb{R}^n) \; \; J_s(u) \stackrel{def}{=} (\mathcal{F}^{-1} \circ A_{\varphi_s} \circ \mathcal{F})(f) = \mathcal{F}^{-1}(\varphi_s \cdot \mathcal{F}(u)).
$$
\end{definition}

Since for any $s \in \mathbb{R}$ the function $\varphi_s$ is an infinitely differentiable function, with this function and all its partial derivatives being of at most polynomial growth, linear operators $\mathbb{A}_{\varphi_s} \colon \mathcal{S}(\mathbb{R}^n) \to \mathcal{S}(\mathbb{R}^n)$ and $A_{\varphi_s} \colon \mathcal{S}'(\mathbb{R}^n) \to \mathcal{S}'(\mathbb{R}^n)$ are well-defined and, as long as for any $s \in \mathbb{R}$ we have
$$
\mathbb{A}_{\varphi_s} \circ \mathbb{A}_{\varphi_{-s}} = \mathbb{A}_{\varphi_{-s}} \circ \mathbb{A}_{\varphi_s} = Id_{\mathcal{S}(\mathbb{R}^n)} \; \; \mbox{and, hence} \; \; A_{\varphi_s} \circ A_{\varphi_{-s}} = A_{\varphi_{-s}} \circ A_{\varphi_s} = Id_{\mathcal{S}'(\mathbb{R}^n)}
$$
it follows that for arbitrary $s \in \mathbb{R}$ the linear operators $\mathbb{J}_s \colon \mathcal{S}(\mathbb{R}^n) \to \mathcal{S}(\mathbb{R}^n)$ and $J_s \colon \mathcal{S}'(\mathbb{R}^n) \to \mathcal{S}'(\mathbb{R}^n)$ are linear homeomorphic automorphisms of, respectively, the Schwartz space $\mathcal{S}(\mathbb{R}^n)$ and the dual Schwartz space $\mathcal{S}'(\mathbb{R}^n)$.

The next theorem will allow us to define a periodic version of the operator $J_s$ for any $s \in \mathbb{R}$.

\begin{theorem}\label{J_s_periodicity}
Let $s \in \mathbb{R}$ and let also $u \in \mathcal{D}'(\mathbb{T}^n)$. 
Then there exists a distribution $v_0 \in \mathcal{S}'_{2 \pi}(\mathbb{R}^n)$ such that the series
$$
\sum\limits_{k \in \mathbb{Z}^n} \left(1 + \left|k\right|^2\right)^{\frac{s}{2}} \cdot c_k(u) \cdot \mathbf{g}_k,
$$
converges to $v_0$ (with respect to the weak--$*$ topology on $\mathcal{S}'(\mathbb{R}^n)$) and
$$
J_s(\widetilde{u}) = v_0.
$$
\end{theorem}

Proof. Since $u \in \mathcal{D}'(\mathbb{T}^n)$, then, by Proposition \ref{periodic_distribution_Fourier_coefficients}, the sequence $(c_k(u) \in \mathbb{C} \: | \; k \in \mathbb{Z}^n)$ has at most polynomial growth and, therefore, the sequence
$$
\left(\left.\left(1 + \left|k\right|^2\right)^{\frac{s}{2}} \cdot c_k(u) \: \right| \; k \in \mathbb{Z}^n\right)
$$
also has at most polynomial growth.

By Corollary \ref{periodic_Schwartz_distribution_Fourier_coefficients}, it follows that there exists a periodic distribution $v_0 \in \mathcal{S}'_{2 \pi}(\mathbb{R}^n)$, such that
$$
v_0 \stackrel{\mathcal{S}'(\mathbb{R}^n)}{=} \sum\limits_{k \in \mathbb{Z}^n} \left(1 + \left|k\right|^2\right)^{\frac{s}{2}} \cdot c_k(u) \cdot \mathbf{g}_k.
$$

So, only the coincidence of the distributions $J_s(\widetilde{u})$ and $v_0$ is left to establish in order to conclude the proof.

By definition of the natural homeomorphism $G \colon \mathcal{D}'(\mathbb{T}^n) \to \mathcal{S}'_{2 \pi}(\mathbb{R}^n)$  we have a validity of the equality
$$
G(u) = \widetilde{u} \stackrel{S'(\mathbb{R}^n)}{=} \sum\limits_{k \in \mathbb{Z}^n} c_k(u) \cdot \mathbf{g}_k
$$
and, because the Fourier transform on $\mathcal{S}'(\mathbb{R}^n)$ is continuous (with respect to its natural weak--$*$ topology), it follows that
$$
\mathcal{F}(\widetilde{u}\,) \stackrel{S'(\mathbb{R}^n)}{=} \sum\limits_{k \in \mathbb{Z}^n} c_k(u) \cdot \mathcal{F}(\mathbf{g}_k).
$$

Let us note that, firstly, the definition of the series $\sum\limits_{k \in \mathbb{Z}^n} c_k(u) \cdot \mathcal{F}(\mathbf{g}_k)$, converging with respect to the weak--$*$ topology on $\mathcal{S}'(\mathbb{R}^n)$, implies that for an arbitrary function $f \in S(\mathbb{R}^n)$ we have
$$
(\varphi_s \cdot \mathcal{F}(\widetilde{u}\,))(f) = (\mathcal{F}(\widetilde{u}\,))(\varphi_s \cdot f) 
= \sum\limits_{k \in \mathbb{Z}^n} c_k(u) \cdot \left(\mathcal{F}(\mathbf{g}_k)\right)(\varphi_s \cdot f),
$$
and, secondly, for any function $\psi \in S(\mathbb{R}^n)$ we have
$$
(\mathcal{F}(\mathbf{g}_k))(\psi) = \mathbf{g}_k(\mathcal{F}(\psi)) = \int\limits_{\mathbb{R}^n} g_k(x) \cdot (\mathcal{F}(\psi))(x) \: d\mu_{cl}(x) =
$$
$$
= \int\limits_{\mathbb{R}^n} \frac{1}{{(2 \pi)^{\frac{n}{2}}}} \cdot e^{i \cdot < k \, ; \: x >} \cdot (\mathcal{F}(\psi))(x) \: d\mu_{cl}(x) = (\mathcal{F}^{-1}(\mathcal{F}(\psi)))(k) = \psi(k) = \delta_{(k)}(\psi),
$$
where $\delta_{(k)} \stackrel{def}{=} T_k(\delta)$ is a shift of the classical delta--function $\delta$ by $k \in \mathbb{Z}^n$ and, as usual in this context, $\delta \in \mathcal{S}'(\mathbb{R}^n)$ is defined by
$$
\forall \: f \in \mathcal{S}(\mathbb{R}^n) \; \; \delta(f) = f(0).
$$

Therefore, for arbitrary function $f \in S(\mathbb{R}^n)$ the validity of the following chain of equalities follows:
$$
(\varphi_s \cdot \mathcal{F}(\widetilde{u} \,))(f) = \sum\limits_{k \in \mathbb{Z}^n} c_k(u) \cdot \left(\mathcal{F}(\mathbf{g}_k)\right)(\varphi_s \cdot f) =
$$
$$
= \sum\limits_{k \in \mathbb{Z}^n} c_k(u) \cdot \varphi_s(k) \cdot f(k) = \sum\limits_{k \in \mathbb{Z}^n} \left(1 + \left|k\right|^2\right)^{\frac{s}{2}} \cdot c_k(u) \cdot f(k).
$$

So let us now fix an arbitrary function $g \in \mathcal{S}(\mathbb{R}^n)$ and prove that
$$
\left(J_s\left(\widetilde{u}\right)\right)(g) = v_0(g).
$$

Taking $\mathcal{F}^{-1}(g) \in \mathcal{S}(\mathbb{R}^n)$ as $f$ in the chain of equalities above, we obtain
$$
(J_s(\widetilde{u}\,))(g) = (\mathcal{F}^{-1}(\varphi_s \cdot \mathcal{F}(\widetilde{u}\,)))(g) = (\varphi_s \cdot \mathcal{F}(\widetilde{u}\,))(\mathcal{F}^{-1}(g)) =
$$
$$
= \sum\limits_{k \in \mathbb{Z}^n} \left(1 + \left|k\right|^2\right)^{\frac{s}{2}} \cdot c_k(u) \cdot (\mathcal{F}^{-1}(g))(k) =
$$
$$
= \sum\limits_{k \in \mathbb{Z}^n} \left(1 + \left|k\right|^2\right)^{\frac{s}{2}} \cdot c_k(u) \cdot \frac{1}{(2 \pi)^{\frac{n}{2}}} \cdot \int\limits_{\mathbb{R}^n} e^{i \cdot < k \, ; \: y >} \cdot g(y) \: d\mu_{cl}(y) =
$$
$$
= \sum\limits_{k \in \mathbb{Z}^n} \left(1 + \left|k\right|^2\right)^{\frac{s}{2}} \cdot c_k(u) \cdot \mathbf{g}_k(g).
$$

Having in mind definition of the weak--$*$ convergence in $\mathcal{S}'(\mathbb{R}^n)$, we deduce that
$$
(J_s(\widetilde{u}\,))(g) = v_0(g), \; \: \mbox{where} \quad
v_0 \stackrel{\mathcal{S}'(\mathbb{R}^n)}{=} \sum\limits_{k \in \mathbb{Z}^n} \left(1 + \left|k\right|^2\right)^{\frac{s}{2}} \cdot c_k(u) \cdot \mathbf{g}_k
$$

Since the function $g \in \mathcal{S}(\mathbb{R}^n)$ was taken arbitrarily, we obtain the equality
$$
J_s(\widetilde{u}\,) \stackrel{\mathcal{S}'(\mathbb{R}^n)}{=} v_0.
$$

This completes the proof of Theorem \ref{J_s_periodicity}.

\medskip

\begin{corollary}\label{J_s_periodicity_on_regular_distributions}
Let $s \in \mathbb{R}$ and $f \in \mathcal{D}(\mathbb{T}^n)$. Then:

$1)$ there holds an equality
$$
\widetilde{\mathbf{f}} \stackrel{\mathcal{S}'(\mathbb{R}^n)}{=} \sum\limits_{k \in \mathbb{Z}^n} C_k(f) \cdot \mathbf{g}_k,
$$

$2)$ there exists a function $g \in \mathcal{D}(\mathbb{T}^n)$, such that the following equalities hold true:
$$
\mathbf{g} \stackrel{\mathcal{D}'(\mathbb{T}^n)}{=} \sum\limits_{k \in \mathbb{Z}^n} (1 + |k|^2)^{\frac{s}{2}} \cdot C_k(f) \cdot \mathbf{f}_k \quad \mbox{and} \quad J_s\left(\widetilde{\mathbf{f}}\right) \stackrel{\mathcal{S}'(\mathbb{R}^n)}{=} \widetilde{\mathbf{g}} \, .
$$
\end{corollary}

Proof. Let us prove the statement of the Part $1)$ first.

Using the result of Remark \ref{relation_between_Fourier_coefficients}, 
we obtain a decomposition
$$
\mathbf{f} \stackrel{\mathcal{D}'(\mathbb{T}^n)}{=} \sum\limits_{k \in \mathbb{Z}^n} C_k(f) \cdot \mathbf{f}_k.
$$

But validity of this decomposition and definition of the linear homeomorphism $G \colon \mathcal{D}'(\mathbb{T}^n) \to \mathcal{S}'_{2 \pi}(\mathbb{R}^n), \quad u \longmapsto \widetilde{u}$, immediately imply that
$$
\widetilde{\mathbf{f}} \stackrel{\mathcal{S}'(\mathbb{R}^n)}{=} \sum\limits_{k \in \mathbb{Z}^n} C_k(f) \cdot \mathbf{g}_k \: .
$$

The proof of Part $1)$ is completed.

Let us now proceed to the proof of Part $2)$. 

Since Remark \ref{Fourier_coefficients_fast_convergence} implies that the sequence $(C_k(f) \: | \;  k \in \mathbb{Z}^n)$, which consists of the Fourier coefficients of $f \in \mathcal{D}(\mathbb{T}^n)$, is a subpolynomially fast decreasing sequence and, hence, the sequence $\left(\left(1 + \left|k\right|^2\right)^{\frac{s}{2}} \cdot C_k(f) \: | \; k \in \mathbb{Z}^n\right)$ also meets the same condition (i.e., it decreases faster than any negative power of $|k|$ as  $|k| \to +\infty$). Therefore, once again applying the result of Remark \ref{Fourier_coefficients_fast_convergence}, there exists a function $g \in \mathcal{D}(\mathbb{T}^n)$ such that its Fourier coefficients satisfy the equality
$$
\forall \: k \in \mathbb{Z}^n \; \; \: C_k(g) = \left(1 + \left|k\right|^2\right)^{\frac{s}{2}} \cdot C_k(f),
$$
and, by Remark \ref{relation_between_Fourier_coefficients}, we obtain the equality
$$
\mathbf{g} \stackrel{\mathcal{D}'(\mathbb{T}^n)}{=} \sum\limits_{k \in \mathbb{Z}^n} \left(1 + \left|k\right|^2\right)^{\frac{s}{2}} \cdot C_k(f) \cdot \mathbf{f}_k \: .
$$

By Theorem \ref{J_s_periodicity}, applied to the distribution $\widetilde{\mathbf{f}} \stackrel{\mathcal{S}'(\mathbb{R}^n)}{=} \sum\limits_{k \in \mathbb{Z}^n} C_k(f) \cdot \mathbf{g_k}$, we arrive at
$$
J_s\left(\,\widetilde{\mathbf{f}}\,\right) \stackrel{\mathcal{S}'(\mathbb{R}^n)}{=} \sum\limits_{k \in \mathbb{Z}^n} \left(1 + \left|k\right|^2\right)^{\frac{s}{2}} \cdot C_k(f) \cdot \mathbf{g}_k \stackrel{\mathcal{S}'(\mathbb{R}^n)}{=} \widetilde{\mathbf{g}}.
$$

The proof of Part $2)$ is completed.

This concludes the proof of Corollary \ref{J_s_periodicity_on_regular_distributions}.

\bigskip

Now we are ready to define a linear operator $J_{s, \: \pi} \colon \mathcal{D}'(\mathbb{T}^n) \to \mathcal{D}'(\mathbb{T}^n)$.

\begin{definition}\label{periodic_J_s_definition}
Let $s \in \mathbb{R}$. Then a linear operator $J_{s, \: \pi} \colon \mathcal{D}'(\mathbb{T}^n) \to \mathcal{D}'(\mathbb{T}^n)$ is defined as follows:
$$
\forall \: u \in \mathcal{D}'(\mathbb{T}^n) \; \; J_{s, \: \pi} (u) \stackrel{\mathcal{D}'(\mathbb{T}^n)}{=} \sum\limits_{k \in \mathbb{Z}^n} \left(1 + \left|k\right|^2\right)^{\frac{s}{2}} \cdot c_k(u) \cdot \mathbf{f}_k.
$$
\end{definition}

Essentially, Theorem \ref{J_s_periodicity} establishes the commutativity of the diagram
$$
\begin{CD}
\mathcal{D}'(\mathbb{T}^n) @>J_{s, \, \pi}>> \mathcal{D}'(\mathbb{T}^n)\\
@VVGV @AAG^{-1}A\\
\mathcal{S}'_{2 \cdot \pi}(\mathbb{R}^n) @>J_s>> \mathcal{S}'_{2 \cdot \pi}(\mathbb{R}^n)
\end{CD}
$$
where $G$ is a natural homeomorphism from $\mathcal{D}'(\mathbb{T}^n)$ onto $\mathcal{S}'_{2 \cdot \pi}(\mathbb{R}^n)$, defined above in terms of the Fourier coefficients, and $G^{-1} \colon \mathcal{S}'_{2 \cdot \pi}(\mathbb{R}^n) \to \mathcal{D}'(\mathbb{T}^n)$ is an inverse linear homeomorphism.

It turns out that the linear operator $J_{s, \: \pi} \colon \mathcal{D}'(\mathbb{T}^n) \to \mathcal{D}'(\mathbb{T}^n)$ has a Schauder basis of its eigenvectors.

\begin{proposition}\label{J_s_pi_eigenfunctions}
Let $s \in \mathbb{R}$. Then the countable family of distributions $(\mathbf{f}_k \: | \; k \in \mathbb{Z}^n)$ forms a Schauder basis of $\mathcal{D}'(\mathbb{T}^n)$ (with respect to its natural weak--$*$ topology), which consists of the eigenfunctions for $J_{s, \: \pi} \colon \mathcal{D}'(\mathbb{T}^n) \to \mathcal{D}'(\mathbb{T}^n)$ as for arbitrary multi--index $m \in \mathbb{Z}^n$ we have
$$
J_{s, \: \pi}(\mathbf{f}_m) = \left(1 + \left|m\right|^2\right)^{\frac{s}{2}} \cdot \mathbf{f}_m.
$$
\end{proposition}

Proof. The fact that the family of regular distributions $(\mathbf{f}_m \: | \; m \in \mathbb{Z}^n)$ forms a Schauder basis  in $\mathcal{D}'(\mathbb{T}^n)$ (with respect to the natural weak--$*$ topology on $\mathcal{D}'(\mathbb{T}^n)$) was already discussed in Remark \ref{D'_Schauder_basis}.

On the other hand, for arbitrary $m \in \mathbb{Z}^n$ validity of an equality
$$
J_{s, \: \pi}(\mathbf{f}_m) = \left(1 + \left|m\right|^2\right)^{\frac{s}{2}} \cdot \mathbf{f}_m
$$
follows immediately from Definition \ref{periodic_J_s_definition} and the fact that for any $m \in \mathbb{Z}^n$ and $k \in \mathbb{Z}^n$ we have the following chain of equalities:
$$
c_k(\mathbf{f}_m) = \mathbf{f}_m(f_{-k}) = \int\limits_{[-\pi \, ; \: \pi]^n} g_m(x) \cdot g_{-k}(x) \: d\mu_{cl}(x) =
$$
$$
= \int\limits_{[-\pi \, ; \: \pi]^n} \frac{1}{(2 \pi)^n} \cdot e^{i \cdot < m - k \, ; \: x >} \: d\mu_{cl}(x) = \begin{cases}
1, \; \; \mbox{если} \: k = m;\\
0, \; \; \mbox{если} \: k \neq m.
\end{cases}
$$

This concludes the proof of Proposition \ref{J_s_pi_eigenfunctions}.

\medskip

Now let us prove a simple proposition, which states the fundamental semigroup condition for the family of linear operators $\{J_{s, \: \pi} \: | \; s \in \mathbb{R}\}$ and some simple implications of this condition.

\begin{proposition}\label{periodic_J_s_basic_properties}
Let $s \in \mathbb{R}$ and $t \in \mathbb{R}$. Then:

$1) \quad J_{s, \, \pi} \circ J_{t, \, \pi} = J_{t, \, \pi} \circ J_{s, \, \pi} = J_{s + t, \: \pi}$;

$2) \quad J_{0, \, \pi} = Id_{\mathcal{D}'(\mathbb{T}^n)}$;

$3)$ the linear operator $J_{-s, \, \pi} \colon \mathcal{D}'(\mathbb{T}^n) \to \mathcal{D}'(\mathbb{T}^n)$ is an inverse to the linear operator $J_{s, \, \pi} \colon \mathcal{D}'(\mathbb{T}^n) \to \mathcal{D}'(\mathbb{T}^n)$.
\end{proposition}

Proof. Let us fix an arbitrary distribution $u \in \mathcal{D}'(\mathbb{T}^n)$.

Then
$$
u \stackrel{\mathcal{D}'(\mathbb{T}^n)}{=} \sum\limits_{k \in \mathbb{Z}^n} c_k(u) \cdot \mathbf{f}_k
$$
and, hence, we have
$$
(J_{s, \, \pi} \circ J_{t, \, \pi})(u) = J_{s, \, \pi}(J_{t, \, \pi}(u)) = J_{s, \, \pi}\left(\sum\limits_{k \in \mathbb{Z}^n} \left(1 + \left|k\right|^2\right)^{\frac{t}{2}} \cdot c_k(u) \cdot \mathbf{f}_k\right) =
$$
$$
= \sum\limits_{k \in \mathbb{Z}^n} \left(1 + \left|k\right|^2\right)^{\frac{s}{2}} \cdot \left(1 + \left|k\right|^2\right)^{\frac{t}{2}} \cdot c_k(u) \cdot \mathbf{f}_k = \sum\limits_{k \in \mathbb{Z}^n} \left(1 + \left|k\right|^2\right)^{\frac{s + t}{2}} \cdot c_k(u) \cdot \mathbf{f}_k = J_{s + t, \, \pi}(u),
$$
and, swapping $s$ and $t$ in this calculation, we also obtain
$$
(J_{t, \, \pi} \circ J_{s, \, \pi})(u) = J_{t + s, \, \pi}(u) = J_{s + t, \, \pi}(u) = (J_{s, \, \pi} \circ J_{t, \, \pi})(u),
$$
which implies the validity of Part $1)$ because of an arbitrary choice of a distribution $u \in \mathcal{D}'(\mathbb{T}^n)$.

The validity of the chain of equalities
$$
J_{0, \, \pi}(u) = \sum\limits_{k \in \mathbb{Z}^n} \left(1 + \left|k\right|^2\right)^0 \cdot c_k(u) \cdot \mathbf{f}_k = \sum\limits_{k \in \mathbb{Z}^n} c_k(u) \cdot \mathbf{f}_k = u
$$
immediately implies the coincidence of the linear operators $J_{0, \, \pi}$ and $Id_{\mathcal{D}'(\mathbb{T}^n)}$, which completes the proof of Part $2)$.

Finally, from the already proved results of Parts $1)$ and $2)$ we obtain the chain of equalities
$$
J_{s, \, \pi} \circ J_{-s, \, \pi} = J_{-s, \, \pi} \circ J_{s, \, \pi} = J_{0, \, \pi} = Id_{\mathcal{D}'(\mathbb{T}^n)},
$$
which means exactly that both operators $J_{s, \, \pi}$ and $J_{-s, \, \pi}$ are not only invertible but they are each other's inverses.

This concludes the proof of Part $3)$ and, hence, of Theorem \ref{periodic_J_s_basic_properties}.

\bigskip

\bigskip

\bigskip

\bigskip

\bigskip

{\Large 5. The scale of the periodic Bessel potential space $H^s_p(\mathbb{T}^n)$: a lifting property of $J_{s, \: \pi}$, duality and embedding theorems.}

\bigskip

\bigskip


Now we introduce the scale of the periodic Bessel potential spaces.

\begin{definition}\label{H^0_p_def}
Let $p \in [1; \: +\infty)$. Then we say that a distribution $u \in \mathcal{D}'(\mathbb{T}^n)$ belongs to the set $H^0_p(\mathbb{T}^n)$, if there exists a function $f \colon \mathbb{T}^n \to \mathbb{C}$, such that $f \in \mathcal{L}_p(\mathbb{T}^n)$ and $u = \mathbf{f}$, where the latter means that
$$
\forall \: \varphi \in \mathcal{D}(\mathbb{T}^n) \; \; \; u(\varphi) = \int\limits_{[-\pi \, ; \: \pi]^n} f_{\pi}(x) \cdot \varphi_{\pi}(x) \: d\mu_{cl}(x).
$$
\end{definition}

\begin{definition}\label{H^s_p_def}
Let $s \in \mathbb{R}$ and $p \in [1; \: +\infty)$. Then
$$
H^s_p(\mathbb{T}^n) \stackrel{def}{=} \left\{ u \in \mathcal{D}'(\mathbb{T}^n) \: | \; J_{s, \, \pi}(u) \in H^0_p(\mathbb{T}^n) \right\}.
$$
\end{definition}

Now for arbitrary $s \in \mathbb{R}$ and $p \in [1; \: +\infty)$ we want to define a norm on the periodic Bessel potential space $H^s_p(\mathbb{T}^n)$, defined above.

Let $s \in \mathbb{R}$ and $p \geqslant 1$. Then we define the function $\| \cdot \|_{H^s_p(\mathbb{T}^n)} \colon H^s_p(\mathbb{T}^n) \to \mathbb{R}$ as follows:

$1)$ in the case $s = 0$ for arbitrary function $f \in \mathcal{L}_p(\mathbb{T}^n)$ we let
$$
\| \mathbf{f} \|_{H^0_p(\mathbb{T}^n)} = \| f \|_{\mathcal{L}_p(\mathbb{T}^n)} = \| [f] \|_{L_p(\mathbb{T}^n)};
$$

$2)$ in the case $s \neq 0$ for arbitrary distribution $u \in H^s_p(\mathbb{T}^n)$ we let
$$
\| u \|_{H^s_p(\mathbb{T}^n)} = \| J_{s, \, \pi}(u) \|_{H^0_p(\mathbb{T}^n)}.
$$

It is easy to check that for arbitrary numbers $s \in \mathbb{R}$ and $p \in [1; \: +\infty)$:

$1)$ the set $H^s_p(\mathbb{T}^n)$ becomes a complex linear space with respect to the naturally defined pointwise linear operations
$$
+ \colon H^s_p(\mathbb{T}^n) \times H^s_p(\mathbb{T}^n) \to H^s_p(\mathbb{T}^n) \quad \mbox{and} \quad \cdot \colon \mathbb{C} \times H^s_p(\mathbb{T}^n) \to H^s_p(\mathbb{T}^n);
$$

$2)$ this complex linear space $(H^s_p(\mathbb{T}^n), \: +, \: \cdot)$ constitutes a normed space with respect to the norm  $\| \cdot \|_{H^s_p(\mathbb{T}^n)} \colon H^s_p(\mathbb{T}^n) \to \mathbb{R}$, defined as follows:

$a)$ in the case $s = 0$ for an arbitrary function $f \in \mathcal{L}_p(\mathbb{T}^n)$ we let
$$
\| \mathbf{f} \|_{H^0_p(\mathbb{T}^n)} = \| f \|_{\mathcal{L}_p(\mathbb{T}^n)};
$$

$b)$ in the case $s \neq 0$ for an arbitrary distribution $u \in H^s_p(\mathbb{T}^n)$ we let
$$
\| u \|_{H^s_p(\mathbb{T}^n)} = \| J_{s, \, \pi}(u) \|_{H^0_p(\mathbb{T}^n)}.
$$

\bigskip

Let us also note that obviously for an arbitrary $p \in [1; \: +\infty)$ a linear mapping $\Gamma_p \colon L_p(\mathbb{T}^n) \to \mathcal{D}'(\mathbb{T}^n)$, which maps an arbitrary equivalence class $[f] \in L_p(\mathbb{T}^n)$ (which implies $f \in \mathcal{L}_p(\mathbb{T}^n)$) to a regular distribution $\mathbf{f} \in \mathcal{D}'(\mathbb{T}^n)$, is an isometric isomorphism between the normed spaces $\left(L_p(\mathbb{T}^n), \: \| \cdot \|_{L_p(\mathbb{T}^n)}\right)$ and $\left(H^0_p(\mathbb{T}^n), \: \| \cdot \|_{H^0_p(\mathbb{T}^n)}\right)$.

Existence of this canonical isometric isomorphism between these two normed spaces allows us to describe the norm $\| \cdot \|_{H^s_2(\mathbb{T}^n)}, \: s \in \mathbb{R}$, in terms of the (generalized) Fourier coefficients.

\begin{remark}\label{toric_H_s_2_coefficients}
Let $s \in \mathbb{R}$ and let also $u \in \mathcal{D}'(\mathbb{T}^n)$. Let us demonstrate that
$$
u \in H^s_2(\mathbb{T}^n) \; \; \mbox{if and only if the series} \; \sum\limits_{k \in \mathbb{Z}^n} \left|c_k(u)\right|^2 \cdot \left(1 + \left|k\right|^2\right)^s \; \mbox{converges},
$$
and for any $u \in H^s_2(\mathbb{T}^n)$ we have
$$
\| u \|_{H^s_2(\mathbb{T}^n)} = \sqrt{\sum\limits_{k \in \mathbb{Z}^n} \left|c_k(u)\right|^2 \cdot \left(1 + \left|k\right|^2\right)^s} \; .
$$

Indeed, let $u \in H^s_2(\mathbb{T}^n)$. Then there exists a function $f \in \mathcal{L}_2(\mathbb{T}^n)$, such that
$$
J_{s, \: \pi}(u) = \mathbf{f}
$$ and, because of Remark \ref{relation_between_Fourier_coefficients}, for any $k \in \mathbb{Z}^n$ we have
$$
C_k(f) = c_k(\mathbf{f}) = c_k(J_{s, \: \pi}(u)) = \left(1 + \left|k\right|^2\right)^{\frac{s}{2}} \cdot c_k(u).
$$
Since for any function $ \varphi \in \mathcal{L}_2(\mathbb{T}^n)$ we have the equality
$$
\left(\| \varphi \|_{\mathcal{L}_2(\mathbb{T}^n)}\right)^2 = \sum\limits_{k \in \mathbb{Z}^n} \left(\left|C_k(\varphi)\right|^2\right)
$$
the series $\sum\limits_{k \in \mathbb{Z}^n} \left(1 + \left|k\right|^2\right)^s \cdot \left|c_k(u)\right|^2$ converges and, moreover, we have
$$
\| u \|_{H^s_2(\mathbb{T}^n)} = \| \mathbf{f} \|_{H^0_2(\mathbb{T}^n)} = \| f \|_{\mathcal{L}_2(\mathbb{T}^n)} = \sqrt{\sum\limits_{k \in \mathbb{Z}^n} \left|C_k(f)\right|^2} = \sqrt{\sum\limits_{k \in \mathbb{Z}^n} \left(1 + \left|k\right|^2\right)^s \cdot \left|c_k(u)\right|^2}.
$$

If, conversely, for some distribution $u \in \mathcal{D}'(\mathbb{T}^n)$ we have the convergence of the series $\sum\limits_{k \in \mathbb{Z}^n} \left(1 + \left|k\right|^2\right)^s \cdot \left|c_k(u)\right|^2$, then, by F. Riesz -- Fischer theorem, the (functional) series
$$
\sum\limits_{k \in \mathbb{Z}^n} \left(1 + \left|k\right|^2\right)^{\frac{s}{2}} \cdot c_k(u) \cdot f_k
$$
converges in $\mathcal{L}_2(\mathbb{T}^n)$ to some function $g \in \mathcal{L}_2(\mathbb{T}^n)$, which, by Remark \ref{relation_between_Fourier_coefficients}, implies that
$$
\mathbf{g} \stackrel{\mathcal{D}'(\mathbb{T}^n)}{=} \sum\limits_{k \in \mathbb{Z}^n} \left(1 + \left|k\right|^2\right)^{\frac{s}{2}} \cdot c_k(u) \cdot \mathbf{f}_k = \sum\limits_{k \in \mathbb{Z}^n} c_k(J_{s, \, \pi}(u)) \cdot \mathbf{f}_k \stackrel{\mathcal{D}'(\mathbb{T}^n)}{=} J_{s, \: \pi}(u)
$$
and, therefore, $u \in H^s_2(\mathbb{R}^n)$.
\end{remark}

The classical result on the lifting property of the operator's $J_{\alpha}$ action on the scale of the Bessel potential spaces $(H^s_p(\mathbb{R}^n) \: | \; s \in \mathbb{N}, \: p \in [1; \: +\infty))$, also can be generalized to the periodic case. Namely, we have the following lemma, which is well--known in the mathematical folklore but we shall give its proof for the sake of exposition's clarity and comprehensiveness.

\begin{lemma}\label{J_s_pi_isometry}
Let $s \in \mathbb{R}$ and $p \geqslant 1$. Let also $\alpha \in \mathbb{R}$. Then the restriction of the linear operator $J_{\alpha, \, \pi} \colon \mathcal{D}'(\mathbb{T}^n) \to \mathcal{D}'(\mathbb{T}^n)$ on $H^s_p(\mathbb{T}^n)$ is an isometric isomorphism between the normed spaces $(H^s_p(\mathbb{T}^n), \: \| \cdot \|_{H^s_p(\mathbb{T}^n)})$ and $(H^{s - \alpha}_p(\mathbb{T}^n), \: \| \cdot \|_{H^{s - \alpha}_p(\mathbb{T}^n)})$.
\end{lemma}

Proof. Let us fix an arbitrary distribution $u_0 \in H^s_p(\mathbb{T}^n)$. Then the distribution $J_{\alpha, \, \pi}(u_0)$ is well--defined as an element of $\mathcal{D}'(\mathbb{T}^n)$, and by the semigroup property of the operator family $\{J_{t, \, \pi} \: | \; t \in \mathbb{R}\}$, which was proven in Part 3) of Proposition \ref{periodic_J_s_basic_properties}, we obtain that
$$
J_{s - \alpha, \, \pi}(J_{\alpha, \, \pi}(u_0)) = J_{s, \, \pi}(u_0) \in H^0_p(\mathbb{T}^n),
$$
which immediately implies the fact that the distribution $J_{\alpha, \, \pi}(u_0)$ is an element of $H^{s - \alpha}_p(\mathbb{T}^n)$. Therefore, the restriction of $J_{\alpha, \, \pi} \colon \mathcal{D}'(\mathbb{T}^n) \to \mathcal{D}'(\mathbb{T}^n)$ on $H^s_p(\mathbb{T}^n)$ is well--defined as a linear operator from $H^s_p(\mathbb{T}^n)$ to $H^{s - \alpha}_p(\mathbb{T}^n)$.

Moreover, using this semigroup property once again, we see that
$$
\| J_{\alpha, \, \pi}(u_0) \|_{H^{s - \alpha}_p(\mathbb{T}^n)} = \| J_{s - \alpha, \, \pi}(J_{\alpha, \, \pi}(u_0)) \|_{H^0_p(\mathbb{T}^n)} = \| J_{s, \, \pi}(u_0) \|_{H^0_p(\mathbb{T}^n)} = \| u_0 \|_{H^s_p(\mathbb{T}^n)},
$$
which implies isometricity of the mapping $J_{\alpha, \, \pi} \colon H^s_p(\mathbb{T}^n) \to H^{s - \alpha}_p(\mathbb{T}^n)$ with respect to the metrics, generated by the norms $\| \cdot \|_{H^s_p(\mathbb{T}^n)}$ and $\| \cdot \|_{H^{s - \alpha}_p(\mathbb{T}^n)}$.

So, only the set--theoretic equality $Im_{J_{\alpha, \, \pi}}(H^s_p(\mathbb{T}^n)) = H^{s - \alpha}_p(\mathbb{T}^n)$ is left to prove.

Let us fix an arbitrary distribution $v_0 \in H^{s - \alpha}_p(\mathbb{T}^n)$. Then, since we already know that $J_{-\alpha, \, \pi}$ maps $H^{s - \alpha}_p(\mathbb{T}^n)$ to $H^s_p(\mathbb{T}^n)$, we deduce that $J_{-\alpha, \, \pi}(v_0) \in H^s_p(\mathbb{T}^n)$. Therefore, by Part $2)$ of Lemma \ref{periodic_J_s_basic_properties}, we obtain
$$
J_{\alpha, \, \pi}(J_{-\alpha, \: \pi}(v_0)) = Id_{\mathcal{D}'(\mathbb{T}^n)}(v_0) = v_0.
$$
Since the choice of $v_0 \in H^{s - \alpha}_p(\mathbb{T}^n)$ was arbitrary, this implies that
$$
H^{s - \alpha}_p(\mathbb{T}^n) \subset Im_{J_{\alpha, \, \pi}}(H^s_p(\mathbb{T}^n)) 
$$
and, as the inverse embedding was also established, we arrive at a set--theoretic equality
$$
Im_{J_{\alpha, \, \pi}}(H^s_p(\mathbb{T}^n)) = H^{s - \alpha}_p(\mathbb{T}^n).
$$

Thus, we proved that the mapping $\left. J_{\alpha, \, \pi} \right|_{H^s_p(\mathbb{T}^n)} \colon H^s_p(\mathbb{T}^n) \to H^{s - \alpha}_p(\mathbb{T}^n)$ is indeed an isometric isomorphism of the corresponding spaces with respect to their natural norms $\| \cdot \|_{H^s_p(\mathbb{T}^n)}$ and $\| \cdot \|_{H^{s - \alpha}_p(\mathbb{T}^n)}$.

This concludes the proof of Lemma \ref{J_s_pi_isometry}.

\bigskip

Similarly to the case of the Bessel potential spaces $(H^s_p(\mathbb{R}^n) \: | \; s \in \mathbb{N}, \: p \in [1; \: +\infty))$, completeness of any normed space $(H^s_p(\mathbb{T}^n), \: \| \cdot \|_{H^s_p(\mathbb{T}^n)})$ (where $s \in \mathbb{R}$ and $p \in [1; \: +\infty)$) with respect to the metric, generated by this norm $\| \cdot \|_{H^s_p(\mathbb{T}^n)}$, immediately follows from the completeness of the normed space $(L_p(\mathbb{T}^n), \: \| \cdot \|_{L_p(\mathbb{T}^n)})$, which is naturally isometrically isomorphic to $(H^0_p(\mathbb{T}^n), \: \| \cdot \|_{H^0_p(\mathbb{T}^n)})$, and the fact that the restriction of $J_{s, \, \pi} \colon H^s_p(\mathbb{T}^n) \to H^0_p(\mathbb{T}^n)$ is an isometric isomorphism with respect to the natural norms of the corresponding spaces, which can be seen as a partial case of Lemma \ref{J_s_pi_isometry}.

Let us show that for arbitrary $s \in \mathbb{R}$ and $p \geqslant 1$ the set of all regular smooth distributions $\mathbf{D}(\mathbb{T}^n)$, which is defined by
$$
\mathbf{D}(\mathbb{T}^n) \stackrel{def}{=} \{ u \in \mathcal{D}'(\mathbb{T}^n) \: | \; \exists \: f \in \mathcal{D}(\mathbb{T}^n) \colon u = \mathbf{f} \},
$$
is dense in $H^s_p(\mathbb{T}^n)$ with respect to the topology, generated by the standard norm $\| \cdot \|_{H^s_p(\mathbb{T}^n)}$ (see also [14, Theorem 3.5.1(ii)]).

First of all, for any $s \in \mathbb{R}$ from the Part $2)$ of Corollary \ref{J_s_periodicity_on_regular_distributions} and the commutative diagram
$$
\begin{CD}
\mathcal{D}'(\mathbb{T}^n) @>J_{-s, \, \pi}>> \mathcal{D}'(\mathbb{T}^n)\\
@VVGV @AAG^{-1}A\\
\mathcal{S}'_{2 \pi}(\mathbb{R}^n) @>J_{-s}>> \mathcal{S}'_{2 \pi}(\mathbb{R}^n)
\end{CD}
$$
we obtain the validity of the embedding
$$
Im_{J_{-s, \, \pi}}(\mathbf{D}(\mathbb{T}^n)) \subset \mathbf{D}(\mathbb{T}^n).
$$

Secondly, as for any $p \in [1; \: +\infty)$ the set $\mathcal{D}(\mathbb{T}^n)$ is dense in $\left(L_p(\mathbb{T}^n), \: \| \cdot \|_{L_p(\mathbb{T}^n)}\right)$, we arrive at a density of $\mathbf{D}(\mathbb{T}^n)$ in $\left(H^0_p(\mathbb{T}^n), \: \| \cdot \|_{H^0_p(\mathbb{T}^n)}\right)$.
 
Therefore, a density of $\mathbf{D}(\mathbb{T}^n)$ in the space $\left(H^s_p(\mathbb{T}^n), \: \| \cdot \|_{H^s_p(\mathbb{T}^n)}\right)$ for arbitrary indices $s \in \mathbb{R}$ and $p \in [1; \: +\infty)$ readily follows from the fact that the restriction of $J_{-s, \, \pi} \colon \mathcal{D}'(\mathbb{T}^n) \to \mathcal{D}'(\mathbb{T}^n)$ on $H^0_p(\mathbb{T}^n)$ is an isometric isomorphism between $H^0_p(\mathbb{T}^n)$ and $H^s_p(\mathbb{T}^n)$ with respect to the norms $\| \cdot \|_{H^0_p(\mathbb{T}^n)}$ and $\| \cdot \|_{H^s_p(\mathbb{T}^n)}$ (which constitutes a partial case of Lemma \ref{J_s_pi_isometry}).

\bigskip

Analogously to the non--periodic situation, for the periodic Bessel potential spaces there can be introduced a notion of duality between $H^s_p(\mathbb{T}^n)$ and $H^{-s}_{p'}(\mathbb{T}^n)$ for any $s \in \mathbb{R}$ and $p \in (1; \: +\infty)$. This construction is well-known but, in order to preserve the unity of notation and also to prove some minor technical assertions let us examine this construction in some detail.

In the sequel for arbitrary $p \in (1; \: +\infty)$ by $p'$ we shall denote its Lebesgue conjugate, i.e. such a number from $(1; \: +\infty)$ that the equality
$$
\frac{1}{p} + \frac{1}{p'} = 1
$$
holds true.

So for arbitrary $s \in \mathbb{R}$ and $p \in (1; \: +\infty)$ we introduce the duality
$$
< \cdot \, ; \: \cdot >_{s, \pi} \colon H^{-s}_{p'}(\mathbb{T}^n) \times H^s_p(\mathbb{T}^n) \to \mathbb{C},
$$
defined by
$$
\forall \: u \in H^{-s}_{p'}(\mathbb{T}^n), \: \forall \: v \in H^s_p(\mathbb{T}^n) \; \; < u \, ; \: v >_{s, \pi} \: = \: < J_{-s, \, \pi}(u), J_{s, \, \pi}(v) >_{0, \pi},
$$
where
$$
\forall \: f \in \mathcal{L}_{p'}(\mathbb{T}^n), \: \forall \: g \in \mathcal{L}_p(\mathbb{T}^n) \; \; < \mathbf{f} \, ; \: \mathbf{g} >_{0, \pi} \: = \: \int\limits_{[-\pi; \: \pi]^n} f_{\pi}(x) \cdot \overline{g_{\pi}(x)} \; d\mu_{cl}(x).
$$

Linearity of the mapping $J_{\alpha} \colon \mathcal{D}'(\mathbb{T}^n) \to \mathcal{D}'(\mathbb{T}^n)$ directly implies the semilinearity of the function $< \cdot \, ; \: \cdot >_{s, \pi} \colon H^{-s}_{p'}(\mathbb{T}^n) \times H^s_p(\mathbb{T}^n) \to \mathbb{C}$ (where $s \in \mathbb{R}$ and $p \in (1; \: +\infty)$):
\begin{gather*}
1) \; \forall \: u_1 \in H^{-s}_{p'}(\mathbb{T}^n), \; \forall \: u_2 \in H^{-s}_{p'}(\mathbb{T}^n), \; \forall \: v \in H^s_p(\mathbb{T}^n) \\
< u_1 + u_2 \, ; \: v >_{s, \pi} \: = \: < u_1 \, ; \: v >_{s, \pi} + < u_2 \, ; \: v >_{s, \pi};\\
2) \; \forall \: u \in H^{-s}_{p'}(\mathbb{T}^n), \; \forall \: v_1 \in H^s_p(\mathbb{T}^n), \; \forall \: v_2 \in H^s_p(\mathbb{T}^n) \\
< u \, ; \: v_1 + v_2 >_{s, \pi} \: = \: < u \, ; \: v_1 >_{s, \pi} + < u \, ; \: v_2 >_{s, \pi};\\
3) \; \forall \: \alpha \in \mathbb{C}, \; \forall \: u \in H^{-s}_{p'}(\mathbb{T}^n), \; \forall \: v \in H^s_p(\mathbb{T}^n) \quad < \alpha \cdot u \, ; \: v >_{s, \pi} \: = \, \alpha \cdot < u \, ; \: v >_{s, \pi};\\
4) \; \forall \: \alpha \in \mathbb{C}, \; \forall \: u \in H^{-s}_{p'}(\mathbb{T}^n), \; \forall \: v \in H^s_p(\mathbb{T}^n) \quad < u \, ; \: \alpha \cdot v >_{s, \pi} \: = \, \overline{\alpha} \, \cdot < u \, ; \: v >_{s, \pi}.
\end{gather*}

\begin{remark}\label{dual_scalar_product_estimate}
Let $s \in \mathbb{R}$ and $p > 1$. Then, from the definition of the duality
$$
< \cdot \, ; \cdot >_{s, \: \pi} \colon H^{-s}_{p'}(\mathbb{T}^n) \times H^s_p(\mathbb{T}^n) \to \mathbb{C}
$$
and the classical H\"older's inequality for the periodic Lebesgue spaces, it follows that
$$
\forall \: u \in H^{-s}_{p'}(\mathbb{T}^n), \; \forall \: v \in H^s_p(\mathbb{T}^n) \quad \; |< u \, ; \: v >_{s, \pi}| \: \leqslant \: \| u \|_{H^{-s}_{p'}(\mathbb{T}^n)} \cdot \| v \|_{H^s_p(\mathbb{T}^n)}.
$$
\end{remark}

Further we need some additional technical results, which are rather trivial but nevertheless will be supplied with complete proofs for the sake of exposition's integrity.

\begin{lemma}\label{zero_dual_scalar_product_representation}
Let $p \in (1; \: +\infty)$. Let also $u \in H^0_p(\mathbb{T}^n)$ and $f \in \mathcal{D}(\mathbb{T}^n)$. Then
$$
1) \; \; < u \, ; \: \mathbf{f} >_{0, \pi} \: = u\left(\overline{f}\right) \quad \mbox{and} \quad \: < \mathbf{f} \, ; \: u >_{0, \pi} \: = \overline{u\left(\overline{f}\right)};
$$
$$
2) \; \; < \mathbf{f} \, ; \: u >_{0, \pi} \: = \: \overline{< u \, ; \: \mathbf{f} >_{0, \pi}}.
$$
\end{lemma}

Proof. Since $u \in H^0_p(\mathbb{T}^n)$, there exists a function $g \in\mathcal{L}_p(\mathbb{T}^n)$, such that $u = \mathbf{g}$. Taking into account the fact that the condition $f \in \mathcal{D}(\mathbb{T}^n)$ implies $\mathbf{f} \in H^0_{p'}(\mathbb{T}^n)$, we immediately arrive at the following chains of equalities:
$$
< u \, ; \: \mathbf{f} >_{0, \pi} \: = \int\limits_{[-\pi; \: \pi]^n} g_{\pi}(x) \cdot \overline{f_{\pi}(x)} \: d\mu_{cl}(x) = \int\limits_{[-\pi; \: \pi]^n} g_{\pi}(x) \cdot \overline{f}_{\pi}(x) \: d\mu_{cl}(x) = u(\overline{f})
$$
\centerline{and}
$$
< \mathbf{f} \, ; \: u >_{0, \pi} \: = \int\limits_{[-\pi; \: \pi]^n} f_{\pi}(x) \cdot \overline{g_{\pi}(x)} \: d\mu_{cl}(x) = \overline{\int\limits_{[-\pi; \: \pi]^n} g_{\pi}(x) \cdot \overline{f_{\pi}(x)} \: d\mu_{cl}(x)} = \overline{u\left(\overline{f}\right)}.
$$

The equality $< \mathbf{f} \, ; \: u >_{0, \pi} \: = \: \overline{< u \, ; \: \mathbf{f} >_{0, \pi}}$ then immediately follows.

This concludes the proof of Lemma \ref{zero_dual_scalar_product_representation}.

\bigskip

\begin{corollary}\label{D'_convergence_implies_DSP_conbvergence}
Let $s \in \mathbb{R}$ and $p \in (1; \: +\infty)$. Let also for arbitrary $n \in \mathbb{Z}_+$ $\; u_n \in \mathcal{D}'(\mathbb{T}^n)$ and, finally, let the sequence of distributions $(u_n \: | \; n \in \mathbb{N})$ converge to a distribution $u_0$ with respect to $\tau_{\mathcal{D}'(\mathbb{T}^n)}$. Then for any function $f \in \mathcal{D}(\mathbb{T}^n)$ we have the following convergencies (with respect to the natural Euclidean topology on $\mathbb{C}$):
$$
< u_n \, ; \: \mathbf{f} >_{0, \pi} \; \xrightarrow[n \to \infty]{} \; < u_0 \, ; \: \mathbf{f} >_{0, \pi} \quad \mbox{and} \quad < \mathbf{f} \, ; \: u_n >_{0, \pi} \; \xrightarrow[n \to \infty]{} \; < \mathbf{f} \, ; \: u_0 >_{0, \pi}.
$$
\end{corollary}

Proof. Fix an arbitrary function $f \in \mathcal{D}(\mathbb{T}^n)$. By Lemma \ref{zero_dual_scalar_product_representation} and the definition of weak--$*$ convergence on $\mathcal{D}'(\mathbb{T}^n)$, we arrive at the following relations:
$$
< u_n \, ; \: \mathbf{f} >_{0, \pi} \: = u_n(\overline{f}) \xrightarrow[n \to \infty]{} u_0(\overline{f}) = \: < u_0 \, ; \: \mathbf{f} >_{0, \pi}
$$
\centerline{and}
$$
< \mathbf{f} \, ; \: u_n >_{0, \pi} \: = \overline{u_n\left(\overline{f}\right)} \xrightarrow[n \to \infty]{} \overline{u_0(\overline{f})} = \: < \mathbf{f} \, ; \: u_0 >_{0, \pi}.
$$

This concludes the proof of Corollary \ref{D'_convergence_implies_DSP_conbvergence}.

\bigskip

\begin{proposition}\label{dual_scalar_product_representation}
Let $s \in \mathbb{R}$ and $p \in (1; \: +\infty)$. Let also $u \in H^s_p(\mathbb{T}^n)$. The for any function $f \in \mathcal{D}(\mathbb{T}^n)$ we have
$$
< u \, ; \: \mathbf{f} >_{-s, \:\pi} \: = \: < \mathbf{f} \, ; \: u >_{s, \:\pi} \: = u\left(\overline{f}\right).
$$
\end{proposition}

Proof. Let us fix an arbitrary function $f \in \mathcal{D}(\mathbb{T}^n)$. Then, by Proposition \ref{periodic_distribution_Fourier_coefficients} and Remark \ref{relation_between_Fourier_coefficients}, it follows that
$$
u \stackrel{\mathcal{D}'(\mathbb{T}^n)}{=} \sum\limits_{k \in \mathbb{Z}^n} c_k(u) \cdot \mathbf{f}_k \quad \mbox{and} \quad \mathbf{f} \stackrel{\mathcal{D}'(\mathbb{T}^n)}{=} \sum\limits_{k \in \mathbb{Z}^n} C_k(f) \cdot \mathbf{f}_k \: ,
$$
which immediately implies the validity of the following representations:
$$
J_{s, \, \pi}(u) \stackrel{\mathcal{D}'(\mathbb{T}^n)}{=} \sum\limits_{k \in \mathbb{Z}^n} \left(1 + \left|k\right|^2\right)^{\frac{s}{2}} \cdot c_k(u) \cdot \mathbf{f}_k \; \: \mbox{and} \; \: J_{-s, \, \pi}(\mathbf{f}) \stackrel{\mathcal{D}'(\mathbb{T}^n)}{=} \sum\limits_{k \in \mathbb{Z}^n} \left(1 + \left|k\right|^2\right)^{-\frac{s}{2}} \cdot C_k(f) \cdot \mathbf{f}_k \: .
$$

Taking into account the fact that $\mathbf{f} \in H^{-s}_{p'}(\mathbb{T}^n)$ and using both semilinearity of $< \cdot \, ; \cdot >_{s, \: \pi}$ and Corollary \ref{D'_convergence_implies_DSP_conbvergence}, we obtain the following chain of equalities:
$$
< u \, ; \: \mathbf{f} >_{-s, \pi} \: = \: < J_{s, \, \pi}(u) \, ; \: J_{-s, \, \pi}(\mathbf{f}) >_{0, \: \pi} \: =
$$
$$
= \: < (\mathcal{D}'(\mathbb{T}^n))-\sum\limits_{k \in \mathbb{Z}^n} \left(1 + \left|k\right|^2\right)^{\frac{s}{2}} \cdot c_k(u) \cdot \mathbf{f}_k \, ; \: J_{-s, \, \pi}(\mathbf{f}) >_{0, \: \pi} =
$$
$$
= \sum\limits_{k \in \mathbb{Z}^n} \left(1 + \left|k\right|^2\right)^{\frac{s}{2}} \cdot c_k(u) \cdot < \mathbf{f}_k \, ; \: J_{-s, \, \pi}(\mathbf{f}) >_{0, \: \pi} \: =
$$
$$
= \sum\limits_{k \in \mathbb{Z}^n} \left(1 + \left|k\right|^2\right)^{\frac{s}{2}} \cdot c_k(u) \cdot < \mathbf{f}_k \, ; \: (\mathcal{D}'(\mathbb{T}^n))-\sum\limits_{l \in \mathbb{Z}^n} \left(1 + \left|l\right|^2\right)^{-\frac{s}{2}} \cdot C_l(f) \cdot \mathbf{f}_l >_{0, \: \pi} \: =
$$
$$
= \sum\limits_{k \in \mathbb{Z}^n} \left(1 + \left|k\right|^2\right)^{\frac{s}{2}} \cdot c_k(u) \cdot \left(\sum\limits_{l \in \mathbb{Z}^n} \left(1 + \left|l\right|^2\right)^{-\frac{s}{2}} \cdot \overline{C_l(f)} \, \cdot < \mathbf{f}_k \, ; \: \mathbf{f}_l >_{0, \pi} \right) = \sum\limits_{k \in \mathbb{Z}^n} c_k(u) \cdot \overline{C_k(f)},
$$
where the last equality follows from the fact that for any multi--indices $k \in \mathbb{Z}^n$ and $l \in \mathbb{Z}^n$ we have
$$
< \mathbf{f}_k \, ; \: \mathbf{f}_l >_{0, \pi} = \frac{1}{(2 \cdot \pi)^n} \cdot \int\limits_{[-\pi \, ; \: \pi]^n} e^{i \cdot < k \, ; \: x >} \cdot \overline{e^{i \cdot < l \, ; \: x >}} \: d\mu_{cl}(x) =
$$
$$
= \frac{1}{(2 \cdot \pi)^n} \cdot \int\limits_{[-\pi \, ; \: \pi]^n} e^{i \cdot < k - l \, ; \: x >} \: d\mu_{cl}(x) =
\begin{cases}
1, \; \: \mbox{if} \: k = l;\\
0, \; \: \mbox{if} \: k \neq l;
\end{cases}
$$

Using the result of Lemma \ref{action_of_distribution_on_function} and the equality $\overline{C_k(f)} = C_{-k}\left(\overline{f}\right)$, valid for any $k \in \mathbb{Z}^n$, it follows that
$$
< u \, ; \: \mathbf{f} >_{-s, \pi} \: = \sum\limits_{k \in \mathbb{Z}^n} c_k(u) \cdot C_{-k}\left(\overline{f}\right) = u\left(\overline{f}\right).
$$

Finally, as long as $J_{-s, \, \pi}(\mathbf{f}) \in \mathbf{D}(\mathbb{T}^n)$ and $J_{s, \, \pi}(u) \in H^0_p(\mathbb{T}^n)$, applying part $3)$ of Lemma \ref{zero_dual_scalar_product_representation} we obtain
$$
< \mathbf{f} \, ; \: u >_{s, \pi} \: = < J_{-s, \, \pi}(\mathbf{f}) \, ; \: J_{s, \, \pi}(u) >_{0, \pi} = \overline{< J_{s, \, \pi}(u) \, ; \: J_{-s, \, \pi}(\mathbf{f}) >_{0, \pi}} = \overline{< u \, ; \: \mathbf{f} >_{-s, \: \pi}}.
$$
and, since we've already proven the equality $< u \, ; \: \mathbf{f} >_{-s, \: \pi} = u\left(\overline{f}\right)$, then the equality
$$
< \mathbf{f} \, ; \: u >_{s, \pi} \: = \overline{u\left(\overline{f}\right)}
$$
holds true.

This concludes the proof of Proposition \ref{dual_scalar_product_representation}.

\bigskip

\begin{remark}\label{dual_DSP_estimate}
Let $s \in \mathbb{R}$ and $p \in (1; \: +\infty)$. Then, by Proposition \ref{dual_scalar_product_representation}, it immediately follows that for an arbitrary distribution $u \in H^{-s}_{p'}(\mathbb{T}^n)$ and an arbitrary function $f \in D(\mathbb{T}^n)$ we have
\begin{equation}\label{eq_DSP_estimate}
|u(f)| = \left| < u \, ; \: \overline{f} >_{s, \pi} \right| \leqslant \| u \|_{H^{-s}_{p'}(\mathbb{T}^n)} \cdot \| \mathbf{f} \|_{H^s_p(\mathbb{T}^n)}.
\end{equation}
Since the set $\mathbf{D}(\mathbb{T}^n)$ is dense in the Banach space $(H^s_p(\mathbb{T}^n), \: \| \cdot \|_{H^s_p(\mathbb{T}^n)})$ validity of this estimate allows us to say that arbitrary distribution $u \in H^{-s}_{p'}(\mathbb{T}^n)$ induces a continuous (with respect to topology, generated by the norm $\| \cdot \|_{H^s_p(\mathbb{T}^n)}$) linear functional on $H^s_p(\mathbb{T}^n)$ with its norm being less or equal than $ \| u \|_{H^{-s}_{p'}(\mathbb{T}^n)}$.
\end{remark}

Let us also note that from this estimate \eqref{eq_DSP_estimate} it follows that for arbitrary $s \in \mathbb{R}$ and $p \in (1; \: +\infty)$ the convergence
$$
u_n \xrightarrow[n \to \infty]{\tau_{H^s_p(\mathbb{T}^n)}} u_0,
$$
where for an arbitrary $n \in \mathbb{Z}_+$ we have $u_n \in H^s_p(\mathbb{T}^n)$, implies a weak convergence in the space $\mathcal{D}'(\mathbb{T}^n)$:
$$
\forall \: f \in D(\mathbb{T}^n) \; \; u_n(f) \xrightarrow[n \to \infty]{} u_0(f).
$$

Remark \ref{dual_DSP_estimate} allows the following well--known inversion (see, e.g., [14, 3.5.6]).

\begin{proposition}\label{dual_inverse_estimate}
Let $s \in \mathbb{R}$ and $p \in (1; \: +\infty)$. Let also $u \in D'(\mathbb{T}^n)$. If there exists a constant $C \in (0; \: +\infty)$, such that
$$
\forall \: f \in D(\mathbb{T}^n) \; \; \; | u(f) | \leqslant \; C \cdot \| \mathbf{f} \|_{H^s_p(\mathbb{T}^n)},
$$
then $u \in H^{-s}_{p'}(\mathbb{T}^n)$ and $\| u \|_{H^{-s}_{p'}(\mathbb{T}^n )} \leqslant C$.
\end{proposition}

In what follows for arbitrary  $s \in \mathbb{R}$ and $p \in [1; \: +\infty)$ let us denote by $(H^s_p(\mathbb{T}^n))^{\overline{*}}$ the linear space of all complex--valued antilinear functionals on $H^s_p(\mathbb{T}^n)$, which are continuous with respect to the norm topology of the space $H^s_p(\mathbb{T}^n)$.

Finally, from Remark \ref{dual_DSP_estimate} and Proposition \ref{dual_inverse_estimate} the main result on the duality of the periodic Bessel potential spaces can be deduced.

\begin{theorem}\label{duality_theorem}
Let $s \in \mathbb{R}$ and $p \in (1; \: +\infty)$. Then there exists a linear isomorphism between the spaces $(H^s_p(\mathbb{T}^n))^{\overline{*}}$ and $H^{-s}_{p'}(\mathbb{T}^n)$, isometric with respect to the norms $\| \cdot \|_{(H^s_p(\mathbb{T}^n))^{\overline{*}}}$ and $\| \cdot \|_{H^{-s}_{p'}(\mathbb{T}^n)}$ and such that for any $w \in (H^s_p(\mathbb{T}^n))^{\overline{*}}$ there exists a distribution $v_w \stackrel{def}{=} F(w) \in H^{-s}_{p'}(\mathbb{T}^n)$, satisfying the condition
$$
\forall \: u \in H^s_p(\mathbb{T}^n) \; \; w(u) = \: < v_w \, ; u >_{s, \: \pi} \: .
$$
\end{theorem}

\bigskip

Of fundamental importance is also the analogue of the classical Sobolev embedding theorem for the periodic Bessel potential spaces (see, e.g., [14, 3.5.1, Remark 4] and [3, Lemma 2.5]).

\begin{theorem}\label{periodic_embedding_theorem}
Let $s \in \mathbb{R}, \: t \in \mathbb{R}, \; p \in (1; \: +\infty)$ and $q \in (1; \: +\infty)$. Let also one of the following conditions hold true:
$$
1) \; p \leqslant q, \: s - \frac{n}{p} \geqslant t - \frac{n}{q} \: ; \quad \quad 2) \; p \geqslant q, \: s \geqslant t.
$$
Then there holds a continuous (with respect to the topologies, induced by the norms $\| \cdot \|_{H^s_p(\mathbb{T}^n)}$ and $\| \cdot \|_{H^t_q(\mathbb{T}^n)}$) embedding
$$
H^s_p(\mathbb{T}^n) \subset H^t_q(\mathbb{T}^n).
$$
\end{theorem}

In the sequel we also need a multiplicative estimate, which is crucial for studying the problem of finding a constructive description of the multipliers acting in the scale of periodic Bessel potential spaces. For the similar estimates in the scale of the (nonperiodic) Bessel potential spaces $(H^s_p(\mathbb{R}^n) \: | \; s \in \mathbb{R}, \: p \in [1; \; +\infty))$ see [1, Lemma 6] and the more general setting in [13, Theorem 8.2.1].

\begin{proposition}\label{main_multiplicative_estimate}
Let $p \in (1; \: +\infty), \: q \in (1; \: +\infty)$ and $s \geqslant t \geqslant 0, \: s > \frac{n}{p}$. Let also either $1) \; p \leqslant q, \; s - \frac{n}{p} \geqslant t - \frac{n}{q}$, or  $2) \; p \geqslant q$. Then there exists a constant $C \in (0; \: +\infty)$, such that for any $f \in \mathcal{D}(\mathbb{T}^n)$ and $g \in \mathcal{D}(\mathbb{T}^n)$ we have the following estimate:
$$
\| f \cdot \mathbf{g} \|_{H^t_q(\mathbb{T}^n)} \leqslant C \cdot \| \mathbf{f} \|_{H^s_p(\mathbb{T}^n)} \cdot \| \mathbf{g} \|_{H^t_q(\mathbb{T}^n)} \; .
$$
\end{proposition}

Proof of this estimate is based upon the equivalent definition of the periodic Bessel potential spaces, which employs the localization technique, generated by some smooth finite partition of unity on the torus $\mathbb{T}^n$ (for discussion of this approach in a general setting of a smooth manifold see [12] while the case of the $n$--dimensional torus is treated in, e.g., [11]). This proof shall be given in another paper by the authors, which investigates the connections between these two approaches to the definition of the periodic Bessel potential spaces, defined on the $n$--dimensional torus.

\bigskip

\bigskip

\bigskip

\bigskip

{\Large 6. Multipliers in the scale of the periodic Bessel potential spaces and their constructive description.}

\bigskip

\bigskip

Now let us define the space of multipliers acting in the scale of the periodic Bessel potential spaces.

\begin{definition}\label{localized_space_def}
 Let $s \in \mathbb{R}$ and $p \in [1; \: +\infty)$. Then we shall say that a distribution $u \in \mathcal{D}'(\mathbb{T}^n)$ belongs to the set $H^s_{p, \: loc}(\mathbb{T}^n)$ if and only if for any function $f \in \mathcal{D}(\mathbb{T}^n)$ we have $f \cdot u \in H^s_p(\mathbb{T}^n)$.
\end{definition}

\begin{definition}\label{multipliers_def}
Let $s \in \mathbb{R}, \: t \in \mathbb{R}, \: p \in [1; \: +\infty), \: q \in [1; \: +\infty)$. Then a distribution $u \in H^t_{q, \: loc}(\mathbb{T}^n)$  is called a multiplier from $H^s_p(\mathbb{R}^n)$ to $H^t_q(\mathbb{R}^n)$ (and we denote it by $u \in M_{\pi}[H^s_p(\mathbb{T}^n) \to H^t_q(\mathbb{T}^n)]$) if  there exists a constant $C \in [0; \: +\infty)$, such that for any function $f \in \mathcal{D}(\mathbb{T}^n)$ we have
$$
\| f \cdot u \|_{H^t_q(\mathbb{T}^n)} \leqslant C \cdot \| \mathbf{f} \|_{H^s_p(\mathbb{T}^n)}.
$$
\end{definition}

Just like in the case of the non--periodic Bessel potential space, it is easy to show that, under the assumptions of Definition \ref{multipliers_def}, the set $M_{\pi}[H^s_p(\mathbb{T}^n) \to H^t_q(\mathbb{T}^n)]$ becomes a linear space with respect to the naturally defined linear operations and, moreover, the function $\| \cdot \|_{M_{\pi}[H^s_p(\mathbb{T}^n) \to H^t_q(\mathbb{T}^n)]} \colon M_{\pi}[H^s_p(\mathbb{T}^n) \to H^t_q(\mathbb{T}^n)] \to \mathbb{R}$, defined as follows
$$
\forall \: u \in M_{\pi}[H^s_p(\mathbb{T}^n) \to H^t_q(\mathbb{T}^n)]
$$
$$
\| u \|_{M_{\pi}[H^s_p(\mathbb{T}^n) \to H^t_q(\mathbb{T}^n)]} = \min \{ C \in [0; \: +\infty)( \: | \; \forall \: f \in \mathcal{D}(\mathbb{T}^n) \; \; \| f \cdot u \|_{H^t_q(\mathbb{T}^n)} \leqslant C \cdot \| \mathbf{f} \|_{H^s_p(\mathbb{T}^n)} \},
$$
is a norm on this linear space.

Since we have already proved that for any $\gamma \in \mathbb{R}$ and $r \in [1; \; +\infty)$ the set $\mathbf{D}(\mathbb{T}^n)$ is dense in $H^{\gamma}_r(\mathbb{T}^n)$ with respect to the natural norm topology of $H^{\gamma}_r(\mathbb{T}^n)$, then, under assumptions of Definition \ref{multipliers_def}, any multiplier $u \in M_{\pi}[H^s_p(\mathbb{T}^n) \to H^t_q(\mathbb{T}^n)]$ induces a generalized multiplication operator $M_u \colon H^s_p(\mathbb{T}^n) \to H^t_q(\mathbb{T}^n)$, which is a linear operator, bounded with respect to the norms $\| \cdot \|_{H^s_p(\mathbb{T}^n)}$ and $\| \cdot \|_{H^t_q(\mathbb{T}^n)}$ and meeting the condition
$$
\forall \; f \in \mathcal{D}(\mathbb{T}^n) \; \; \: M_u(\mathbf{f}) = f \cdot u.
$$
Moreover, for any multiplier $u \in M_{\pi}[H^s_p(\mathbb{T}^n) \to H^t_q(\mathbb{T}^n)]$ we have
$$
\| u \|_{M_{\pi}[H^s_p(\mathbb{T}^n) \to H^t_q(\mathbb{T}^n)]} = \| M_u \|_{\mathcal{B}(H^s_p(\mathbb{R}^n), \: H^t_q(\mathbb{R}^n))},
$$
where $\| \cdot \|_{\mathcal{B}(H^s_p(\mathbb{R}^n), \: H^t_q(\mathbb{R}^n))}$ is a natural operator norm, whose value is defined for any linear operator from $H^s_p(\mathbb{T}^n)$ to $H^t_q(\mathbb{T}^n)$ which is bounded with respect to the norms $\| \cdot \|_{H^s_p(\mathbb{R}^n)}$ and $\| \cdot \|_{H^t_q(\mathbb{R}^n)}$.

\medskip

The following result on the interrelation between multipliers and duality also has an analogue in the non-periodic case (see [1, Proposition 1].

\begin{theorem}\label{multipliers_indices_symmetry}
Let $s \in \mathbb{R}, \; t \in \mathbb{R}, \; p \in (1; \: +\infty), \; q \in (1; \: +\infty)$. Then the multiplier spaces $M_{\pi}[H^s_p(\mathbb{T}^n) \to H^t_q(\mathbb{T}^n)]$ and $M_{\pi}[H^{-t}_{q'}(\mathbb{T}^n) \to H^{-s}_{p'}(\mathbb{T}^n)]$ coincide and their standard norms are equal.
\end{theorem}

Proof. Suppose $u \in M_{\pi}[H^s_p(\mathbb{T}^n) \to H^t_q(\mathbb{T}^n)]$.

Now let us fix an arbitrary function $f \in \mathcal{D}(\mathbb{T}^n)$. Since the distribution $u$ is a multiplier from $H^s_p(\mathbb{T}^n)$ to $H^t_q(\mathbb{T}^n)$, then $f \cdot u \in H^t_q(\mathbb{T}^n)$ and also the following estimate holds:
$$
\| f \cdot u \|_{H^t_q(\mathbb{T}^n)} \leqslant \| u \|_{M_{\pi}[H^s_p(\mathbb{T}^n) \to H^t_q(\mathbb{T}^n)]} \cdot \| \mathbf{f} \|_{H^s_p(\mathbb{T}^n)}.
$$
Applying Remark \ref{dual_DSP_estimate}, it follows that
$$
\forall \: g \in \mathcal{D}(\mathbb{T}^n) \; \; \; \left|(g \cdot u)(f)\right| = \left|(f \cdot u)(g)\right| \leqslant \| f \cdot u \|_{H^t_q(\mathbb{T}^n)} \cdot \| \mathbf{g} \|_{H^{-t}_{q'}(\mathbb{T}^n)} \leqslant
$$
$$
\leqslant \| u \|_{M_{\pi}[H^s_p(\mathbb{T}^n) \to H^t_q(\mathbb{T}^n)]} \cdot \| \mathbf{g} \|_{H^{-t}_{q'}(\mathbb{T}^n)} \cdot \| \mathbf{f} \|_{H^s_p(\mathbb{T}^n)}.
$$
Since the choice of the function $f \in \mathcal{D}(\mathbb{T}^n)$ was arbitrary, then, by Proposition \ref{dual_inverse_estimate}, we deduce that for any function $g \in \mathcal{D}(\mathbb{T}^n)$ we have $g \cdot u \in H^{-s}_{p'}(\mathbb{T}^n)$ and
$$
\| g \cdot u \|_{H^{-s}_{p'}(\mathbb{T}^n)} \leqslant \| u \|_{M_{\pi}[H^s_p(\mathbb{T}^n) \to H^t_q(\mathbb{T}^n)]} \cdot \| \mathbf{g} \|_{H^{-t}_{q'}(\mathbb{T}^n)}.
$$
But, by Definition \ref{multipliers_def}, it means exactly that $u \in M_{\pi}[H^{-t}_{q'}(\mathbb{T}^n) \to H^{-s}_{p'}(\mathbb{T}^n)]$ and the following inequality is valid:
$$
\| u \|_{M_{\pi}[H^{-t}_{q'}(\mathbb{T}^n) \to H^{-s}_{p'}(\mathbb{T}^n)]} \leqslant \| u \|_{M_{\pi}[H^s_p(\mathbb{T}^n) \to H^t_q(\mathbb{T}^n)]}.
$$

Therefore, since $u \in M_{\pi}[H^s_p(\mathbb{T}^n) \to H^t_q(\mathbb{T}^n)]$ was taken arbitrarily, the embedding
$$
M_{\pi}[H^s_p(\mathbb{T}^n) \to H^t_q(\mathbb{T}^n)] \subset M_{\pi}[H^{-t}_{q'}(\mathbb{T}^n) \to H^{-s}_{p'}(\mathbb{T}^n)]
$$
is valid and the following estimate holds true:
$$
\forall \: u \in M_{\pi}[H^s_p(\mathbb{T}^n) \to H^t_q(\mathbb{T}^n)] \; \; \; \| u \|_{M_{\pi}[H^{-t}_{q'}(\mathbb{T}^n) \to H^{-s}_{p'}(\mathbb{T}^n)]} \leqslant \| u \|_{M_{\pi}[H^s_p(\mathbb{T}^n) \to H^t_q(\mathbb{T}^n)]}.
$$

Since $-(-s) = s, -(-t) = t, (p')' = p$ and $(q')' = q$, swapping $s$ with $t$ and $p$ with $q$ in the considerations above, we obtain the validity of both the embedding
$$
M_{\pi}[H^{-t}_{q'}(\mathbb{T}^n) \to H^{-s}_{p'}(\mathbb{T}^n)] \subset M_{\pi}[H^s_p(\mathbb{T}^n) \to H^t_q(\mathbb{T}^n)]
$$
and the estimate
$$
\forall \: u \in M_{\pi}[H^{-t}_{q'}(\mathbb{T}^n) \to H^{-s}_{p'}(\mathbb{T}^n)] \; \; \; \| u \|_{M_{\pi}[H^s_p(\mathbb{T}^n) \to H^t_q(\mathbb{T}^n)]} \leqslant \| u \|_{M_{\pi}[H^{-t}_{q'}(\mathbb{T}^n) \to H^{-s}_{p'}(\mathbb{T}^n)]}.
$$

This concludes the proof of Theorem \ref{multipliers_indices_symmetry}.

\bigskip

Now we shall prove our main result which establishes a constructive description of the multiplier space in the case when the so--called Strichartz type conditions on the smoothness indices do hold.

\smallskip

\begin{theorem}\label{periodic_multipliers_description}
Let $s \in [0; \: +\infty), \: t \in [0; \: +\infty), \: p \in (1; \; +\infty), \: q \in (1; \: +\infty)$. Additionally,

if $s \geqslant t$, then let $s > \frac{n}{p}$ and one of the following conditions hold
$$
1) \; p \leqslant q', \; s - \frac{n}{p} \geqslant t - \frac{n}{q'} \quad \mbox{or}  \quad 2) \; p \geqslant q',
$$

and if $t \geqslant s$, then let $t > \frac{n}{q'}$ and one of the following conditions hold
$$
3) \; q' \leqslant p, \; t - \frac{n}{q'} \geqslant s - \frac{n}{p} \quad \mbox{or} \quad 4) \; q' \geqslant p.
$$

Then
$$
M_{\pi}[H^s_p(\mathbb{T}^n) \to H^{-t}_q(\mathbb{T}^n)] = H^{-t}_q(\mathbb{T}^n) \cap H^{-s}_{p'}(\mathbb{T}^n),
$$
and the norms $\| \cdot \|_{M_{\pi}[H^s_p(\mathbb{T}^n) \to H^{-t}_q(\mathbb{T}^n)]}$ and $\| \cdot \|_{H^{-t}_q(\mathbb{T}^n) \cap H^{-s}_{p'}(\mathbb{T}^n)}$ are equivalent. Here we define the norm $\| \cdot \|_{H^{-t}_q(\mathbb{T}^n) \cap H^{-s}_{p'}(\mathbb{T}^n)}$ by letting
$$
\forall \: u \in H^{-t}_q(\mathbb{T}^n) \cap H^{-s}_{p'}(\mathbb{T}^n) \quad \; \| u \|_{H^{-t}_q(\mathbb{T}^n) \cap H^{-s}_{p'}(\mathbb{T}^n)} = \max\left(\| u \|_{H^{-t}_q(\mathbb{T}^n)}; \: \| \cdot \|_{H^{-s}_{p'}(\mathbb{T}^n)}\right).
$$
\end{theorem}

Proof. Let us at first consider the case when $s \geqslant t, \; s > \frac{n}{p}$ and one of the following two conditions holds:
$$
1) \; p \leqslant q', \; s - \frac{n}{p} \geqslant t - \frac{n}{q'}; \; \; \; \; 2) \; p \geqslant q'.
$$

Then, since  $\frac{1}{p} + \frac{1}{p'} = 1 = \frac{1}{q} + \frac{1}{q'}$ implies $\frac{1}{p'} - \frac{1}{q} = \frac{1}{q'} - \frac{1}{p}$, it follows that either
$$
q \leqslant p' \quad \mbox{and} \quad -t - \frac{n}{q} \geqslant -s - \frac{n}{p'}
$$
\centerline{or}
$$
-t \geqslant -s \quad \mbox{and} \quad q \geqslant p'.
$$
Then, by Theorem \ref{periodic_embedding_theorem}, we have a continuous embedding $H^{-t}_q(\mathbb{T}^n) \subset H^{-s}_{p'}(\mathbb{T}^n)$ and, therefore, there exits a constant $K \in [0; \: +\infty)$, such that
$$
\forall \: u \in H^{-t}_q(\mathbb{T}^n) \; \; \| u \|_{H^{-s}_{p'}(\mathbb{T}^n)} \leqslant C \cdot \| u \|_{H^{-t}_q(\mathbb{T}^n)}.
$$
This implies a set--theoretic equality
$$
H^{-t}_q(\mathbb{T}^n) \cap H^{-s}_{p'}(\mathbb{T}^n) = H^{-t}_q(\mathbb{T}^n),
$$
and an equivalence of the norms $\| \cdot \|_{H^{-t}_q(\mathbb{T}^n) \cap H^{-s}_{p'}(\mathbb{T}^n)}$ and $\| \cdot \|_{H^{-t}_q(\mathbb{T}^n)}$.

Now suppose $u \in H^{-t}_q(\mathbb{T}^n)$.

Let us also fix an arbitrary function $f \in \mathcal{D}(\mathbb{T}^n)$. Then for an arbitrary function $g \in \mathcal{D}(\mathbb{T}^n)$ we have $f \cdot g \in \mathcal{D}(\mathbb{T}^n)$ and also, by Remark \ref{dual_DSP_estimate}, we have
$$
\left|(f \cdot u)(g)\right| = \left|u(f \cdot g)\right| \leqslant \| u \|_{H^{-t}_q(\mathbb{T}^n)} \cdot \| f \cdot \mathbf{g} \|_{H^t_{q'}(\mathbb{T}^n)}.
$$
Employing a multiplicative estimate from Proposition \ref{main_multiplicative_estimate}, then we obtain the following inequality:
$$
\forall \: g \in \mathcal{D}(\mathbb{T}^n) \; \; \:  \left|(f \cdot u)(g)\right| \leqslant C_0 \cdot \| u \|_{H^{-t}_q(\mathbb{T}^n)} \cdot \| \mathbf{f} \|_{H^s_p(\mathbb{T}^n)} \cdot \| \mathbf{g} \|_{H^t_{q'}(\mathbb{T}^n)}.
$$
Therefore, by Proposition \ref{dual_inverse_estimate}, we see that $f \cdot u \in H^{-t}_q(\mathbb{T}^n)$ and
$$
\| f \cdot u \|_{H^{-t}_q(\mathbb{T}^n)} \leqslant C_0 \cdot \| u \|_{H^{-t}_q(\mathbb{T}^n)} \cdot \| \mathbf{f} \|_{H^s_p(\mathbb{T}^n)}.
$$
As the function $f \in \mathcal{D}(\mathbb{T}^n)$ was taken arbitrarily, it follows that $u$ belongs to the multiplier space $M_{\pi}[H^s_p(\mathbb{T}^n) \to H^{-t}_q(\mathbb{T}^n)]$ and the estimate
$$
\| u \|_{M_{\pi}[H^s_p(\mathbb{T}^n) \to H^{-t}_q(\mathbb{T}^n)]} \leqslant C_0 \cdot \| u \|_{H^{-t}_q(\mathbb{T}^n)}
$$
holds true.

Because of the arbitrary choice of $u \in H^{-t}_q(\mathbb{T}^n)$, it follows that there holds a continuous (with respect to the corresponding norm topologies) embedding
$$
H^{-t}_q(\mathbb{T}^n) \subset M_{\pi}[H^s_p(\mathbb{T}^n) \to H^{-t}_q(\mathbb{T}^n)].
$$

Let us prove the inverse embedding
$$
M_{\pi}[H^s_p(\mathbb{T}^n) \to H^{-t}_q(\mathbb{T}^n)] \subset H^{-t}_q(\mathbb{T}^n).
$$

Assume $v \in M_{\pi}[H^s_p(\mathbb{T}^n) \to H^{-t}_q(\mathbb{T}^n)]$. 

Since for the function $E \in \mathcal{D}(\mathbb{T}^n)$, defined by
$$
\forall \: z \in \mathbb{T}^n \; \; E(z) = 1,
$$
we have $\mathbf{E} \in H^s_p(\mathbb{T}^n)$, then $M_v(\mathbf{E}) \in H^{-t}_q(\mathbb{T}^n)$.

But taking into account the fact that, as we noted above, $M_v(\mathbf{E}) = E \cdot v = v$, it follows that $v \in H^{-t}_q(\mathbb{T}^n)$ and
$$
\| v \|_{H^{-t}_q(\mathbb{T}^n)} = \| M_v(\mathbf{E}) \|_{H^{-t}_q(\mathbb{T}^n)} \leqslant \| M_v \|_{\mathcal{B}(H^s_p(\mathbb{T}^n), \: H^{-t}_q(\mathbb{T}^n))} \cdot \| \mathbf{E} \|_{H^s_p(\mathbb{T}^n)} =
$$
$$
= \| \mathbf{E} \|_{H^s_p(\mathbb{T}^n)} \cdot \| v \|_{M_{\pi}[H^s_p(\mathbb{T}^n) \to H^{-t}_q(\mathbb{T}^n)]}.
$$

So, we have obtained an estimate 
$$
\| v \|_{H^{-t}_q(\mathbb{T}^n)}\leqslant \| \mathbf{E} \|_{H^s_p(\mathbb{T}^n)} \cdot \| v \|_{M_{\pi}[H^s_p(\mathbb{T}^n) \to H^{-t}_q(\mathbb{T}^n)]},
$$
which, by virtue of an arbitrary choice of $v \in M_{\pi}[H^s_p(\mathbb{T}^n) \to H^{-t}_q(\mathbb{T}^n)]$, means exactly the validity and the continuity (with respect to the corresponding norm topologies) of the embedding
$$
M_{\pi}[H^s_p(\mathbb{T}^n) \to H^{-t}_q(\mathbb{T}^n)] \subset H^{-t}_q(\mathbb{T}^n),
$$
where the needed estimate is valid for the positive constant $\| \mathbf{E} \|_{H^s_p(\mathbb{T}^n)}$.

Thus, we have proved that 
$$
M_{\pi}[H^s_p(\mathbb{T}^n) \to H^{-t}_q(\mathbb{T}^n)] = H^{-t}_q(\mathbb{T}^n)
$$
and the natural norms of these spaces are equivalent.

Let us consider the case $t \geqslant s, \; t > \frac{n}{q'}$ and one of the following conditions holds: $3) \; q' \leqslant p, \; t - \frac{n}{q'} \geqslant s - \frac{n}{p}$ and $4) \; q' \geqslant p$.

Using analogous reasoning as in the previous case, we obtain that one of the following conditions is valid: either $p' \leqslant q$ and $-s - \frac{n}{p'} \geqslant -t - \frac{n}{q}$, or $p' \geqslant q$ and $-s \geqslant -t$. Therefore, the continuous (with respect to the norm topologies of these spaces) embedding $H^{-s}_{p'}(\mathbb{T}^n) \subset H^{-t}_q(\mathbb{T}^n)$ is valid, and it follows that
$$
H^{-s}_{p'}(\mathbb{T}^n) \cap H^{-t}_q(\mathbb{T}^n) = H^{-s}_{p'}(\mathbb{T}^n),
$$
with the norms $\| \cdot \|_{H^{-t}_q(\mathbb{T}^n) \cap H^{-s}_{p'}(\mathbb{T}^n)}$ and $\| \cdot \|_{H^{-s}_{p'}(\mathbb{T}^n)}$ being equivalent.

According to Theorem \ref{multipliers_indices_symmetry}, we arrive at
$$
M_{\pi}[H^s_p(\mathbb{T}^n) \to H^{-t}_q(\mathbb{T}^n)] = M_{\pi}[H^t_{q'}(\mathbb{T}^n) \to H^{-s}_{p'}(\mathbb{T}^n)],
$$
with norms of these spaces being equivalent. Since in our case $t \geqslant s, \: t > \frac{n}{q'}$ and one of the following conditions is valid: $3) \; q' \leqslant p, \; t - \frac{n}{q'} \geqslant s - \frac{n}{p}$ and $4) \; q' \geqslant p$, then, by the result already proved in Part $1)$, it follows that
$$
M_{\pi}[H^t_{q'}(\mathbb{T}^n) \to H^{-s}_{p'}(\mathbb{T}^n)] = H^{-s}_{p'}(\mathbb{T}^n)
$$
and the norms $\| \cdot \|_{M_{\pi}[H^t_{q'}(\mathbb{T}^n) \to H^{-s}_{p'}(\mathbb{T}^n)]}$ and $\| \cdot \|_{H^{-s}_{p'}(\mathbb{T}^n)}$ are equivalent.

Finally, we arrive at a conclusion that
$$
M_{\pi}[H^s_p(\mathbb{T}^n) \to H^{-t}_q(\mathbb{T}^n)] = H^{-s}_{p'}(\mathbb{T}^n)
$$
and the norms $\| \cdot \|_{M_{\pi}[H^s_p(\mathbb{T}^n) \to H^{-t}_q(\mathbb{T}^n)]}$ and $\| \cdot \|_{H^{-s}_{p'}(\mathbb{T}^n)}$ are equivalent.

This completes the proof of Theorem \ref{periodic_multipliers_description}.

\bigskip

\bigskip

{\itshape {\bfseries References.}}

\bigskip

1. Belyaev~A.\,А., Shkalikov~А.\,А., {\itshape Multipliers in spaces of Bessel potentials: the case of indices of nonnegative smoothness,}
Math. Notes, {\bfseries 102} (2017), №\,5, 632 --- 644

\smallskip

2. Belyaev~A.\,А., Shkalikov~А.\,А., {\itshape Multipliers in Bessel potential spaces with \\ smoothness indices of different sign,} St. Petersburg Math. J., {\bfseries 30} (2019), №\,2, 203 --- 218

\smallskip

3. Cirant~M., Goffi~A., {\itshape On the existence and uniqueness of solutions to time-dependent fractional MFG,} SIAM J. Math. Anal., {\bfseries 51} (2019), №\,2, 913 --- 954

\smallskip

4. Colombo~M., de Lellis~C., de Rosa~L., {\itshape Ill-posedness of Leray solutions for the hypodissipative Navier--Stokes equations,} Comm. Math. Phys., {\bfseries 362} (2017), №\,2, 659 --- 688

\smallskip

5. Grafakos~L., {\itshape Classical Fourier Analysis,} New York:  Springer, 2008

\smallskip

6. Grubb~G., {\itshape Distributions and Operators,} New York: Springer, 2009

\smallskip

7. Grubb~G., {\itshape Fractional Laplacians on domains, a development of H\"ormander's theory of $\mu$--transmission pseudodifferential operators,} Adv. Math., {\bfseries 268} (2015), 478 --- 528

\smallskip

8. Haroske~D., Triebel~H., {\itshape Distributions, Sobolev spaces, elliptic equations,} Z\"urich: European Mathematical Society Publishing House, 2007

\smallskip

9. Kappeler~T., M\"ohr~C., {\itshape Estimates for periodic and Dirichlet eigenvalues of the Schroedinger operator with singular potentials,} J. Funct. Anal., {\bfseries 186} (2001), №\,1, 62 --- 91

\smallskip

10. Roncal~L., Stinga~P.\,R.,  {\itshape Transference of fractional Laplacian regularity,} Special Functions, Partial Differential Equations, and Harmonic Analysis, New York: \\ Springer, 2014 (pp. 203 -- 212)

\smallskip

11. Roncal~L., Stinga~P.R., {\itshape Fractional Laplacian on the torus,} Commun. Contemp. Math., {\bfseries 18} (2016), №\,2, 1550033 (26 pages)

\smallskip

12. Rosenberg~S., {\itshape The Laplacian on a Riemannian manifold: An introduction to analysis on manifolds,} Cambridge: Cambridge University Press, 1997

\smallskip

13. Runst~T., Sickel~W., {\itshape Sobolev spaces of fractional order, Nemytskij operators, and nonlinear partial differential equations,} Berlin: De Gruyter, 2011

\smallskip

14. Schmeisser~H.J., Triebel~H., {\itshape Topics in Fourier Analysis and Function Spaces,} Chichester, New York: Wiley, 1987




\bigskip

\bigskip

Addresses:

Department of Mechanics and Mathematics,
Lomonosov Moscow State University,
Moscow, Russia

and

S.\,M.~Nikol'skii Mathematical Institute,
Peoples' Friendship University of Russia,
Moscow, Russia

{email: belyaev\_aa@pfur.ru}

\newpage

\bigskip

\bigskip

Department of Mechanics and Mathematics,
Lomonosov Moscow State University,
Moscow, Russia

{email: shkalikov@mi.ras.ru}

\end{document}